\documentclass[12pt]{article}

\usepackage{amssymb,amsmath,graphicx}
\usepackage{amsfonts,amsthm,color,url}
\usepackage[caption=false]{subfig}
\usepackage{epsfig,latexsym,graphicx}
\usepackage[authoryear]{natbib}
\newtheorem{theorem}{Theorem}[section]

\newtheorem{lemma}[theorem]{Lemma}

\usepackage{setspace}

\begin{document}

\title{\bf{Robust Sure Independence Screening  \\for Non-polynomial dimensional\\ Generalized Linear Models
}}

\author{Abhik Ghosh$^1$, Erica Ponzi$^2$, Torkjel Sandanger$^3$ and  Magne Thoresen$^2*$\\
$^1$ Indian Statistical Institute, Kolkata, India\\
$^2$ University of Oslo, Oslo, Norway\\
$^3$ 
UiT,  The Arctic University of Norway, Troms\o, Norway.\\
*Corresponding Author. Email: magne.thoresen@medisin.uio.no}
\date{}

\maketitle
\noindent
\textbf{Running headline:} Robust sure independence screening for GLMs

\newpage
\begin{abstract}
	We consider the problem of variable screening in ultra-high  dimensional generalized linear models (GLMs) of non-polynomial orders.
Since the popular SIS approach is extremely unstable in the presence of contamination and noise, 
we discuss a new robust screening procedure based on the minimum density power divergence estimator (MDPDE) of the marginal regression coefficients.
Our proposed screening procedure performs well  under pure and contaminated data scenarios. 
We provide a theoretical motivation for the use of marginal MDPDEs for variable screening from both population as well as sample aspects;
in particular, we prove that the marginal MDPDEs are uniformly consistent leading to the sure screening property of our proposed algorithm.
Finally, we propose an appropriate MDPDE based extension for robust conditional screening in GLMs
along with the derivation of its sure screening property. 
Our proposed methods are illustrated through extensive numerical studies along with an interesting real data application. 
\end{abstract}


\section{Introduction}
\label{SEC:intro}

The class of generalized linear models (GLMs) is a rich class of  parametric regression models 
that allows to study a wide range of relationship structures for different types of response data,
which makes the GLMs one of the most popular statistical tools for real-life applications across many disciplines. 
Let us consider the GLM in its canonical form:  
given a set of $p$ predictor variables $X_1, X_2, \ldots, X_p$, 
the scalar response variable $Y$ follows  a distribution from the exponential family having density
\begin{equation}
	f(y;\theta) = \exp\left\{ y\theta - b(\theta) + c(y) \right\},
	\label{EQ:exp_family_density}
\end{equation}
for some appropriate (known) functions $b(\cdot)$ and $c(\cdot)$ and the unknown canonical parameter $\theta$.
For simplicity, we do not consider a dispersion parameter in the model (e.g., logistic or Poison regression) 
although it can easily be incorporated in all our methodological discussions and the theories derived throughout the paper with slight modifications. 
We concentrate on the mean regression model for $\theta$ given by 
\begin{eqnarray}
	E[Y|\boldsymbol{X}=\boldsymbol{x}] = b'(\theta) =  g^{-1}\left(\boldsymbol{x}^T\boldsymbol{\beta}\right),
	\label{EQ:GLM}
\end{eqnarray}
where $\boldsymbol{X}=(X_0=1, X_1, \ldots, X_p)^T$, 
$\boldsymbol{\beta}=(\beta_0, \beta_1, \ldots, \beta_p)^T\in\mathbb{R}^{p+1}$ is the vector of unknown regression coefficients 
and $g$ is a monotone differentiable link function.
Given independent and identically distributed (IID) data $(y_i, \boldsymbol{x}_i)$, $i=1, \ldots, n$, 
our objective is  to fit a GLM by efficiently estimating $\boldsymbol{\beta}$ and use it for subsequent inference.

Commonly, the regression coefficient $\boldsymbol{\beta}$ is estimated through likelihood based approaches (or 
suitable extensions)  under the classical low dimensional set-up ({ $p<n$}).
However, recent advancements of technologies across disciplines generates data on a large number of possible covariates with limited observations
leading to $p\gg n$, known as the high-dimensional set-up. In this paper, we consider ultra-high dimensional  GLMs 
with the number of covariates being of non-polynomial (NP) order of $n$, i.e., $\log(p) = O(n^l) $ for some $0<l<1$.
However, to perform meaningful inference in such situations, we need to assume sparsity of the model --- 
only $s\ll n$ covariates (out of the vast pool of $p$ covariates) are actually important to explain the variability in the response. 
There are several statistical procedures like LASSO or other regularized approaches 
\citep{Fan/Li:2001,Buhlmann/VanDeGeer:2011,Hastie/etc:2015,Giraud:2014,Ghosh/Majumdar:2019} 
to simultaneously select these important variables and estimate the corresponding (non-zero) regression coefficients. 
Although they often work reasonably well in moderately high dimensions, their computation becomes highly extensive in ultra-high dimensional set-ups.
Therefore, it is more efficient to first reduce the set of all covariates to a sufficiently small size (maybe $<n$)
through some initial screening procedure. Among these, the most popular is the sure independence screening (SIS) 
proposed by \cite{Fan/Lv:2008} for the linear regression model and  
later extended to GLMs by  \cite{Fan/Song:2010}. The SIS has become extremely popular for its simplicity, elegance,
computational speed as well as the theoretical guarantees for sure screening of the true model, asymptotically with probability tending to one. 
Subsequently, SIS has been extended to different types of data and associated statistical problems;
see, e.g., \cite{Barut/etc:2016,Zhao/Li:2012,Luo/etc:2014,Saldana/Feng:2018} among many others.

The SIS is commonly applied, besides other applications, in the context of omics data which generally include different types of noise and outliers;
the same issue of data contamination also often arises in other real-life applications involving  extremely large number of features. 
However, the SIS procedure and its extensions are mainly based on the Pearson correlation or the maximum likelihood estimator (MLE) 
of the marginal regression coefficients, both of which are non-robust against possible outliers in the data.
This non-robustness of SIS was, in fact, first noted in the discussion of the original paper itself by \cite{Gather/Guddat:2008}.
They proposed an alternative robust SIS using the Gnanadesikan-Kettenring correlation in place of the usual correlation 
while ranking the covariates in a linear regression model. Subsequently, several other robust versions of SIS,
mostly non-parametric in nature, were proposed for  the high-dimensional linear regression model only
\citep{Hall/Miller:2009,Li/etc:2012a,Li/etc:2012b,Mu/Xiong:2014,Wang/etc:2017}.
Although these non-parametric versions of SIS can potentially be applied to the GLMs as well (possibly with appropriate modulation),
they were never theoretically studied in the literature. 
Thus, there is a need for a robust variable screening procedure for the ultra-high dimensional GLM 
with proper theoretical guarantees of its sure screening property. 
We aim to fill this gap in the literature by developing 
a robust sure screening procedure for the general class of GLMs.

Compared to any non-parametric robust procedure, a parametric robust approach is known to provide significantly 
higher efficiency when the assumed model is valid for a majority of the data except for the noise/contamination part
\citep[e.g.,][]{Hampel/etc:1986}.
Recently, a robust parametric version of SIS, namely the DPD-SIS, 
has been proposed for ultra-high dimensional linear regression models by \cite{Ghosh/Thoresen:2020}.
This DPD-SIS is empirically studied and found to have significantly improved performance compared to the other existing non-parametric SIS procedures 
under data contamination, although no theoretical guarantees are provided in \cite{Ghosh/Thoresen:2020}.
They have proposed to use the marginal regression approach as in \cite{Fan/Song:2010} but to estimate the marginal regression slopes  
by the robust minimum density power divergence estimator (MDPDE) instead of the MLE. 
These MDPDEs were first proposed by \cite{Basu/etc:1998} as a robust generalization of the MLE for simple IID problems.
Due to their high robustness along with their high efficiency and simple computation, 
the MDPDEs are subsequently extended to more complex statistical models.
For linear regression models, the MDPDEs are studied by, e.g., \cite{Ghosh/Basu:2013}. 
For different GLMs as well, the MDPDEs are seen to provide highly efficient and robust parameter estimates
under the classical low dimensional set-ups \citep{Basu/etc:2011,Ghosh/Basu:2016,Basu/etc:2017,Basu/etc:2018,Ghosh:2019}. 
In this paper, we utilize the MDPDEs under the marginal regression approach to develop a robust variable screening procedure 
for ultra-high dimensional GLMs, as an extension of the robust DPD-SIS of \cite{Ghosh/Thoresen:2020}.
Additionally, we prove that the proposed procedure satisfies the sure screening property for the general class of GLMs
and that it can also control the selection of false positives under appropriate assumptions,
which needed quite non-trivial extensions of the existing theories. 
{ To our knowledge, this is the first robust variable screening procedure for the general class of GLMs 
	(beyond simple linear regression) with proper theoretical guarantees.
}

Further, we also extend our proposed DPD-SIS to develop a robust conditional screening procedure
under NP-dimensional GLMs, which we will refer to as the conditional DPD-SIS.
The conditional SIS (CSIS), proposed by \cite{Barut/etc:2016}, has been a natural extension of the usual SIS 
that can take care of additional information (whenever available) about some previously chosen important variables. 
Among several advantages, most importantly, CSIS helps to select the hidden important variables. 
However, just like usual SIS, the CSIS is also extremely non-robust under data contamination, 
and there is no literature available on its (parametric) robust version. Our proposed conditional DPD-SIS serves this purpose. 
Its population-level justifications as well as the sample-level sure screening property 
are also derived rigorously  under reasonably practical assumptions.

We illustrate the proposed DPD-SIS through appropriate simulation studies of ultra-high dimensional GLMs, 
in addition to an interesting real data application involving the search for biomarkers in lung cancer.
For simplicity in presentation, all proofs are deferred to the Appendixes.


\section{The Proposed DPD-SIS for NP-dimensional GLMs}
\label{SEC:DPD_SIS}

Let us consider the GLM described in (\ref{EQ:exp_family_density})-(\ref{EQ:GLM}) with 
ultra-high dimensional covariates;
for simplicity in presentation, throughout the rest of the paper we will assume canonical link function so that $b'= g^{-1}$
and hence $\theta = \boldsymbol{x}^T\boldsymbol{\beta}$ in (\ref{EQ:GLM}).
Suppose that the true value of the regression coefficient $\boldsymbol{\beta}$ is denoted by 
$\boldsymbol{\beta}_0=(\beta_{00}, \beta_{01}, \ldots, \beta_{0p})^T$.
We assume that the true model, denoted as $\mathcal{M}_0 = \left\{ 1\leq j \leq p : \beta_{0j} \neq 0 \right\}$, 
is sparse with model size $s=|\mathcal{M}_0| < n$.
Our aim is to perform an initial screening of the covariates in a robust manner
such that it includes all the truly important variables corresponding to $\mathcal{M}_0$;
this property is referred to as the sure screening property in the literature.

We follow the marginal regression approach of \cite{Fan/Song:2010}
to consider the GLM  for $Y$ based on $X_j$ (plus an intercept term) separately for each $j=1, \ldots, p$;
let us denote the associated regression coefficients for these marginal models by $\boldsymbol{\beta}_j^M = (\beta_{j0}^M, \beta_j^M)$, respectively.
In this context, we also assume that the covariates are standardized 
so that $E(X_j) = 0$ and $E(X_j^2)=1$ for all $j=1, \ldots, p$.
However, instead of using the MLE of $\boldsymbol{\beta}_j^M$ as in \cite{Fan/Song:2010}, 
we propose to use their MDPDEs. 
Note that each marginal GLM is of low dimension,  having only two parameters in $\boldsymbol{\beta}_j^M$. 
Hence, we can follow \cite{Ghosh/Basu:2016} to define  their MDPDE as the minimizer of an appropriately defined average DPD measure between 
the observed data $(y_i, \boldsymbol{x}_i)$, $i=1, \ldots, n$, and the assumed GLM density (\ref{EQ:exp_family_density}).
After simplifications, the MDPDE of $\boldsymbol{\beta}_j^M$ with (given) tuning parameter $\alpha> 0$ is defined as
\begin{eqnarray}
	~~~~~~\widehat{\boldsymbol{\beta}}_j^{M\alpha} = \left(\widehat{\beta}_{j0}^{M\alpha}, \widehat{\beta}_j^{M\alpha}\right )
	=\arg\min\limits_{\beta_{j0}, \beta_j} \frac{1}{n}\sum_{i=1}^n l_\alpha\left(y_i, \beta_{j0} + \beta_j x_{ij}\right),
	\label{EQ:MMDPDE}
\end{eqnarray}
where $\boldsymbol{x}_i = (1, x_{i1}, \ldots, x_{ip})^T$ for each $i=1, \ldots, n$, and 
$l_\alpha(y, \theta) =  \int f(s;\theta)^{1+\alpha} ds - \left(1+\frac{1}{\alpha}\right) f(y;\theta)^\alpha +\frac{1}{\alpha}$. 
The tuning parameter $\alpha$ in the definition of the MDPDE is known to control the trade-off between 
the efficiency under pure data and the robustness against data contamination.
In fact, $l_0(y, \theta) := \lim\limits_{\alpha\rightarrow 0} l_\alpha(y, \theta)
= - \log f(y, \theta)$ so that the MDPDE at $\alpha=0$ (in a limiting sense) is nothing but the most efficient and highly non-robust MLE.
For $\alpha>0$, the MDPDE provides a robust extension of the MLE having increasing robustness with a slight loss in efficiency 
as $\alpha$ increases \citep{Ghosh/Basu:2016}.
For any given $\alpha\geq 0$, each MDPDE $\widehat{\boldsymbol{\beta}}_j^{M\alpha}$
can also be obtained by solving the corresponding estimating equation
$
\sum\limits_{i=1}^n \psi_\alpha(y_i, \beta_0 + \beta_j x_{ij})[1, x_{ij}]^T = \boldsymbol{0}_2,
$
and 
\begin{eqnarray}
	&&\psi_\alpha(y, \theta)=  (y - b'(\theta))f(y; \theta)^{\alpha} - \xi_\alpha(\theta), 
	\label{EQ:DPD_psi}
\end{eqnarray}
where $\xi_\alpha(\theta) = \int (y-b'(\theta))f(y; \theta)^{1+\alpha} dy$. 
Note that $\xi_0(\theta)=0$ and hence $\psi_0(y, \theta) = (y - b'(\theta))$, the usual score function, which again leads to the MLE. 
The above form of $\psi_\alpha$, used in the estimating equation of the MDPDE, indeed justifies the particular choice of this DPD loss function $l_\alpha$
for achieving robust solutions. Note that, the first part of the $\psi_\alpha$ function 
(except $\xi_\alpha$, which is there to make the estimating equation unbiased at the model) is basically a weighted score function
with the weight being $f(y; \theta)^{\alpha}$. So, for any $\alpha>0$, observations in the region of low model probability 
(the outlying observations coming from data contamination) get downweighted with regard to their contributions to the MDPDE estimating equation;
as a consequence, the effects of such outlying observations get reduced leading to robust parameter estimates and subsequently to robust screening results
{ (see Section \ref{SEC:IF} for further discussions on robustness
and Section \ref{SEC:tuning} for comments on the choice of  $\alpha$).
}

Based on the MDPDEs $\widehat{\boldsymbol{\beta}}_j^{M\alpha}$ for the marginal regression coefficients, for each $j=1, \ldots, p$,
and a given $\alpha>0$, we choose a suitable pre-defined threshold $\gamma_n$ and select the variables in the set 
\begin{eqnarray}
	\widehat{M}_\alpha(\gamma_n) = \left\{ 1\leq j \leq p : \left|\widehat{\beta}_j^{M\alpha}\right| \geq \gamma_n  \right\}.
	\label{EQ:DPD_SIS_set}
\end{eqnarray}  
By choosing $\gamma_n$ appropriately, we can reduce the number of covariates from a large $p$ to 
any smaller target, say $d<n$, so that the subsequent computation becomes feasible.  
With these $d$  variables selected in $\widehat{M}_\alpha(\gamma_n)$, 
we can then fit any appropriate low-dimensional estimation procedure or regularized approach to get 
the estimated coefficient vector, say 
$\widehat{\boldsymbol{\beta}}_d =(\widehat{\beta}_{d0}, \widehat{\beta}_{d1}, \ldots, \widehat{\beta}_{dd})^T$
and subsequently the final model $\widehat{\mathcal{M}} =  \left\{ 1\leq j \leq p : \widehat{\beta}_{dj} \neq 0 \right\}$.

%
%
%

In our DPD-SIS, we suggest to choose $\gamma_n$ from the target of attaining a fixed model size. 
However, in practice, it may be chosen in several other ways, e.g., controlling the false positives, or prediction error, etc. 
(see Section \ref{SEC:prop_DPDSIS} for an optimal rate of $\gamma_n$).
Note that the case $\alpha=0$  reduces to the ordinary SIS.
Along the same lines, we will show the sure screening property of our proposed DPD-SIS, 
so that  asymptotically   $\widehat{M}_\alpha(\gamma_n)$ contains the true model $\mathcal{M}_0$ with probability tending to one, 
for any given $\alpha>0$. 


\section{Sure Screening property of the DPD-SIS} 
\label{SEC:prop_DPDSIS}

We are considering the GLM in (\ref{EQ:exp_family_density})-(\ref{EQ:GLM}) with the canonical link function 
and the true sparse parameter value $\boldsymbol{\beta}_0$ having support $\mathcal{M}_0$ of  size $s=|\mathcal{M}_0| < n$, 
as described in previous sections.
Recall that, under the ultra-high dimensional set-up considered here, the number of covariates $p=p_n$ is assumed to grow exponentially 
with the sample size $n$; we also allow the true model size $s=s_n$ to depend on $n$. 
Further we assume that the data $(y_i, \boldsymbol{x}_i)$, for $i=1, \ldots, n$, are  IID from a true joint distribution 
$\Pi(dy, d\boldsymbol{x}) = F_{\boldsymbol{\beta}_0}(dy|\boldsymbol{x}) Q(d\boldsymbol{x})$,
where $F_{\boldsymbol{\beta}_0}$ is the (conditional) distribution corresponding to the GLM in (\ref{EQ:exp_family_density})--(\ref{EQ:GLM}) 
and $Q$ is the marginal distribution of the covariates (for which no model is assumed).
Then, it is straightforward from the definition of $\psi_\alpha$ in (\ref{EQ:DPD_psi}) that, for any $\alpha\geq 0$,  
\begin{eqnarray}
	E[\psi_\alpha(Y,\boldsymbol{x}^T\boldsymbol{\beta}_0)|\boldsymbol{X}=\boldsymbol{x}] 
	=E[\psi_\alpha(Y,\boldsymbol{x}_1^T\boldsymbol{\beta}_{01})|\boldsymbol{X}=\boldsymbol{x}] = \boldsymbol{0},
	\label{EQ:Est_eqn_Mpop0}
\end{eqnarray}
where we have used the notation $\boldsymbol{x}_1 = (x_j : j\in \mathcal{M}_0)^T$ and 
$\boldsymbol{\beta}_{01} = (\beta_{0j} : j\in \mathcal{M}_0)^T$. 
Without loss of generality, let us assume $\mathcal{M}_0 = \{1, 2, \ldots, s\}$ 
and consider the partitions $\boldsymbol{x}^T=(\boldsymbol{x}_1^T, \boldsymbol{x}_2^T)$,
$\boldsymbol{\beta}_0^T=(\boldsymbol{\beta}_{01}^T, \boldsymbol{\beta}_{02}^T)$, 
$\boldsymbol{x}_i^T=(\boldsymbol{x}_{i1}^T, \boldsymbol{x}_{i2}^T)$, 
$\boldsymbol{\beta}^T = (\boldsymbol{\beta}_1^T, \boldsymbol{\beta}_2^T)$ and so on for any $p$-vectors,
where the first partitions (e.g., $\boldsymbol{x}_1$, $\boldsymbol{\beta}_{01}$, $\boldsymbol{x}_{i1}$,  $\boldsymbol{\beta}_1$, etc.) 
are of length $s$.

\subsection{Population level results}
\label{SEC:prop_DPDSIS_pop}

We first investigate the proposed DPD-SIS at population level.
The population version (functional) of the marginal MDPDE $\widehat{\boldsymbol{\beta}}_j^{M\alpha}$, defined in (\ref{EQ:MMDPDE}),
is given by  
\begin{eqnarray}
	\boldsymbol{\beta}_j^{M\alpha} = ({\beta}_{j0}^{M\alpha}, \beta_j^{M\alpha})
	=\arg\min\limits_{\beta_{j0}, \beta_j} E \left[l_\alpha\left(Y, \beta_{j0} + \beta_j X_{j}\right)\right].
\end{eqnarray}
This marginal MDPDE functional $\boldsymbol{\beta}_j^{M\alpha}$ then satisfies the estimating equations 
\begin{eqnarray}
	E \left[\psi_\alpha\left(Y, \beta_0 + \beta_j X_{j}\right)\right] = 0,
	~~~~~~~~
	E \left[\psi_\alpha\left(Y, \beta_0 + \beta_j X_{j}\right)X_j\right] = 0.
	\label{EQ:Est_eqn_Mpop}
\end{eqnarray}
Let us denote
$B_\alpha(v(\boldsymbol{x}))= b'(\boldsymbol{x}^T\boldsymbol{\beta}_0) - E[\psi_\alpha(Y, v(\boldsymbol{x}))|\boldsymbol{X}=\boldsymbol{x}].$
Clearly $B_\alpha(\boldsymbol{x}^T\boldsymbol{\beta}_0) = b'(\boldsymbol{x}^T\boldsymbol{\beta}_0)$
but $B_\alpha(\beta_{j0}^{M\alpha} + \beta_j^{M\alpha}x_j)$ do not necessarily equal $b'(\beta_{j0}^{M\alpha} + \beta_j^{M\alpha}x_j)$.
However, at $\alpha=0$ we always have $B_0(v(\boldsymbol{x}))= b'(v(\boldsymbol{x}))$ for any $v(\boldsymbol{x})$.
Using equations (\ref{EQ:Est_eqn_Mpop0}) and (\ref{EQ:Est_eqn_Mpop}), we have proved the following two theorems.
They show why the proposed DPD-SIS is expected to have the targeted sure screening property, at population level;
the proofs are given in Appendix B for brevity. 

\begin{theorem}\label{THM:DPD_SIS_Pop1}
	For a given $\alpha\geq 0$, 
	and for any $j=1, \ldots, p$, the marginal MDPDE functional $\beta_j^{M\alpha}=0$ 
	if and only if Cov$(b'(\boldsymbol{X}^T\boldsymbol{\beta}_0), X_j)=\mbox{Cov}(Y, X_j)=0$.
\end{theorem}

\bigskip
\begin{theorem}\label{THM:DPD_SIS_Pop2}
	Given any $\alpha\geq 0$, { and $j\in \mathcal{M}_0$, 
	assume that either (B1) or (B2)  holds: 
	\begin{itemize}
		\item[(B1)]  $B_\alpha'(\cdot)$ is bounded.
		\item[(B2)] $B_\alpha(t)$ is strictly increasing in $t$ and  $G_\alpha(|x|) = \sup_{|u|\leq |x|}|B_\alpha(u)|$ satisfies 
		\begin{eqnarray}
			E[G_{\alpha}(a|X_j|)|X_j|I(|X_j|\geq n^\eta)] \leq dn^{-\kappa},
			~\mbox{ for some constants } a, d>0, ~\eta\in (0,\kappa). ~~
			\label{EQ:cond_prop2.1}
		\end{eqnarray}
	\end{itemize}
Then, whenever }there exists a constant $c_1>0$ such that 
	$\left|\mbox{Cov}(b'(\boldsymbol{X}^T\boldsymbol{\beta}_0), X_j)\right|\geq c_1n^{-\kappa}$, 
	we have 
	$
	\min_{j\in\mathcal{M}_0}|\beta_j^{M\alpha}| \geq c_2 n^{-\kappa}, 
	$ 
	for some constant $c_2>0$.
\end{theorem}

The above two theorems are similar to Theorems 2 and 3 of \cite{Fan/Song:2010} from the context of SIS,
although the assumptions in our Theorem \ref{THM:DPD_SIS_Pop2} are required on the quantity $B_\alpha(\cdot)$  instead of $b'(\cdot)$
which coincide at $\alpha=0$.  
For any $\alpha\geq 0$, one can indeed show that $B_\alpha'(\cdot)$ is bounded for the normal and the logistic regression models
{ whereas Condition (\ref{EQ:cond_prop2.1}) holds for Poisson regression with suitable covariates; see Appendix \ref{APP:Assumption}}. 

In the same spirit of \cite{Fan/Song:2010}, Theorem \ref{THM:DPD_SIS_Pop1} implies that if 
the set of unimportant covariates $\{X_j : j \notin\mathcal{M}_0\}$ is independent of the  set of important covariates $\{X_j : j\in \mathcal{M}_0\}$
then $\beta_j^{M\alpha} = 0$ for all $j \notin \mathcal{M}_0$ and all $\alpha\geq 0$. 
Further, note that, an important covariate $X_j$ having non-zero correlation with the response has 
a marginal regression coefficient $\beta_j^{M\alpha}\neq 0$.  
These together indicate the existence of a threshold $\gamma_n$ satisfying $\min_{j\in\mathcal{M}_0} |\beta_j^{M\alpha}| \geq \gamma_n$ 
and $\max_{j\notin \mathcal{M}_0} |\beta_j^{M\alpha}|=0$. This forms the theoretical basis for the model selection consistency of the
proposed DPD-SIS with any $\alpha\geq 0$ and justifies our proposal as an authenticate screening criterion.

On the other hand, Theorem \ref{THM:DPD_SIS_Pop2} provides the conditions to yield
$\min_{j\in\mathcal{M}_0} |\beta_j^{M\alpha}| \geq O(n^{-\kappa})$ for some $\kappa<1/2$, 
which can be interpreted as the marginal signals being stronger than the maximum stochastic noise level.
It is an intuitive necessity for the proposed DPD-SIS, the sample version (\ref{EQ:DPD_SIS_set}),
to have the sure screening property. Other than the $\alpha$-dependent assumption,
one crucial condition in Theorem \ref{THM:DPD_SIS_Pop2} is 
$\left|\mbox{Cov}(b'(\boldsymbol{X}^T\boldsymbol{\beta}_0), X_j)\right|\geq c_1n^{-\kappa}$ 
which is the same as required by the usual SIS in \cite{Fan/Song:2010}; 
it can be further simplified for jointly Gaussian covariates following Proposition 1 of \cite{Fan/Song:2010} and the discussion thereafter. 
This theorem also provides the necessary framework to achieve sparsity in the final selected model (\ref{EQ:DPD_SIS_set}).

\subsection{Sample level results}
\label{SEC:prop_DPDSIS_sample}


We first show that the marginal MDPDEs $\widehat{\boldsymbol{\beta}}_j^{M\alpha}$, $j=1, \ldots, p$, are uniformly consistent at an exponential rate
which leads to the sure screening property (sample level) of our proposed DPD-SIS. 
In this regard, let us note that the marginal MDPDE functional $\boldsymbol{\beta}_j^{M\alpha}$ is unique and is an interior point of the parameter space
by convexity of the DPD loss function $l_\alpha\left(Y, \beta_0 + \beta_j X_{j}\right)$ in $\boldsymbol{\beta}_j=(\beta_{0j}, \beta_j)$ for each $j$.
So, we can restrict the minimization of the marginal DPD loss function over the compact set $\mathcal{B}=\{|\beta_{j0}|\leq B, |\beta_j|\leq B \}$
for some large enough constant $B>0$ such that $\boldsymbol{\beta}_j^{M\alpha}$ is also  an interior point of $\mathcal{B}$. 
For each $j=1, \ldots, p$, let $\boldsymbol{X}_j=(1, X_j)^T$ and define the matrices \\
\begin{eqnarray}
	\boldsymbol{J}_{j,\alpha}(\boldsymbol{\beta}_j) &=& E\left[\nabla^2 l_\alpha\left(Y, \beta_0 + \beta_j X_{j}\right)\right]
	= (1+\alpha)E\left[\Gamma_\alpha\left(\boldsymbol{X}_j^T\boldsymbol{\beta}_j\right)\boldsymbol{X}_j\boldsymbol{X}_j^T\right], 
	\label{EQ:J}\\
	\boldsymbol{K}_{j,\alpha}(\boldsymbol{\beta}_j) 
	&=& E\left[\left(\nabla l_\alpha\left(Y, \beta_0 + \beta_j X_{j}\right)\right)\left(\nabla l_\alpha\left(Y, \beta_0 + \beta_j X_{j}\right)\right)^T\right] 
	\nonumber\\
	&=& (1+\alpha)^2E\left[
	\left\{\Gamma_{2\alpha}\left(\boldsymbol{X}_j^T\boldsymbol{\beta}_j\right) - \xi_\alpha^2\left(\boldsymbol{X}_j^T\boldsymbol{\beta}_j\right)\right\}
	\boldsymbol{X}_j\boldsymbol{X}_j^T\right],
	\label{EQ:K}
\end{eqnarray}
where $\Gamma_\alpha(\theta) = \int (y-b'(\theta))^2f(y; \theta)^{1+\alpha} dy$.
Then, the following assumptions are needed for our subsequent theoretical investigation of the DPD-SIS;
here $\alpha\geq 0$ is a fixed given tuning parameter
and $\Lambda_{\min}[\cdot]$ and $\Lambda_{\max}[\cdot]$, respectively, denote the minimum 
and maximum eigenvalues of its argument matrix.

\begin{enumerate}
	\item[(A1)]  The GLM is such that the density $f^{\alpha}$ in (\ref{EQ:exp_family_density}) is bounded { by some constant $L_\alpha>0$,} 
	and $b''(\cdot)$ is continuous and positive. Also, { $|\xi_{\alpha}(\theta)|$ is non-decreasing in $\theta$}.
	
	\item[(A2)] For all $\boldsymbol{\beta}_j\in\mathcal{B}$, there exists some constant $V>0$ such that 
	$\Lambda_{\min}\left[\boldsymbol{J}_{j,\alpha}(\boldsymbol{\beta}_j)\right]\geq V $ uniformly over $j=1, \ldots, p$.   
	
	\item[(A3)]  $\boldsymbol{K}_{j,\alpha}(\boldsymbol{\beta}_j^{M\alpha})$ is finite and positive definite for each $j=1, \ldots, p$.
	Also, the norm $	||\boldsymbol{K}_{j,\alpha}(\boldsymbol{\beta}_j)||_\mathcal{B} = \sup\limits_{\boldsymbol{\beta}_j\in\mathcal{B}, ||\boldsymbol{u}||=1}
	||\boldsymbol{K}_{j,\alpha}(\boldsymbol{\beta}_j)^{1/2}\boldsymbol{u}||$ is bounded from above for each $j$.
	
	\item[(A4)]  There exist an $\epsilon_1>0$ and a large constant $K_n>0$, such that 
	$$
	\sup\limits_{\boldsymbol{\beta}_j\in\mathcal{B}: ||\boldsymbol{\beta}_j-\boldsymbol{\beta}_j^{M\alpha}||\leq \epsilon_1}
	E\left[|B_\alpha(\boldsymbol{X}_j^T\boldsymbol{\beta}_j)|||\boldsymbol{X}_j||_2I(|X_j|>K_n)\right]\leq o\left(\frac{1}{n}\right), 
	~~~~~\mbox{ for all } j=1, 2, \ldots, p.
	$$
	
	\item[(A5)] The distribution of the covariate $X_j$ is such that, for sufficiently large $t>0$ and some positive constants $m_0, m_1, m_2, m_3$ and $\tau$,
	we have
	$
	P(|X_j|>t) = (m_1-m_2)e^{-m_0t^\tau},
	$
	for $j=1, 2, \ldots, p$, and
	$$
	E\left[ \exp\left(b(\boldsymbol{X}^T\boldsymbol{\beta}_0+m_3) - b(\boldsymbol{X}^T\boldsymbol{\beta}_0)\right)\right]
	+ 	E\left[ \exp\left(b(\boldsymbol{X}^T\boldsymbol{\beta}_0-m_3) - b(\boldsymbol{X}^T\boldsymbol{\beta}_0)\right)\right] \leq m_2.
	$$
	
	\item[(A6)] $Var(\boldsymbol{X}^T\boldsymbol{\beta}_0)$ is bounded both from below and above
	by finite positive constants.
	
	\item[(A7)] Either $b''(\cdot)$ is bounded or $\widetilde{\boldsymbol{X}}=(X_1, \ldots, X_p)^T$ 
	follows an elliptically contoured distribution with variance $\Sigma_1$ and
	$\left|E\left[b'(\boldsymbol{X}^T\boldsymbol{\beta}_0)(\boldsymbol{X}^T\boldsymbol{\beta}_0 - \beta_{00})\right]\right|$ is bounded.	 
\end{enumerate}

Note that Assumptions (A1)--(A6) are appropriate extensions of the assumptions made by \cite{Fan/Song:2010} 
to prove the sure screening property of the usual SIS; they coincide at $\alpha=0$ since $L_0=1$, $\xi_0 \equiv 0$, 
$\boldsymbol{J}_{j,0}=\boldsymbol{K}_{j,0}=E\left[b''\left(\boldsymbol{X}_j^T\boldsymbol{\beta}_j\right)\boldsymbol{X}_j\boldsymbol{X}_j^T\right]$
and $B_0(v(\boldsymbol{x}))= b'(v(\boldsymbol{x}))$ for any $v(\boldsymbol{x})$. 
For any $\alpha>0$, Assumption (A1) clearly holds for most common examples of GLM including the normal, Poisson and logistic regression models;
other assumptions are also valid for these GLMs under mild sufficient conditions.
Interestingly, Assumptions (A5)--(A7) are independent of the choice of $\alpha$ and are exactly the same as Assumptions (D), (F) and (G) 
of \cite{Fan/Song:2010}, respectively. In particular, Assumption (A5) ensures that the covariates and the response variable have light tails;  
it implies, via Lemma 1 of \cite{Fan/Song:2010}, that
 
\begin{eqnarray}
	P\left(|Y|\geq \frac{m_0}{m_3}t^{\tau}\right) \leq m_2 e^{-m_0t^\tau}, ~~~~\mbox{for any }~t>0.
	\label{EQ:Tail_response}
\end{eqnarray} 
Assumption (A6), on the other hand, implies that the variance of the response $Y$ is bounded. 
In fact, denoting the variance of $\boldsymbol{X}$ by $\Sigma=\mbox{Diag}\{0, \Sigma_1\}$, Assumption (A6) states that  
$
Var(\boldsymbol{X}^T\boldsymbol{\beta}_0)=\boldsymbol{\beta}_0^T\Sigma\boldsymbol{\beta}_0 =O(1).
$
Note that the maximum eigenvalue of $\Sigma_1$ in Assumption (A7) is the same as $\Lambda_{\max}(\Sigma)$, 
which is a positive finite number by Assumption (A6). 
Further, Assumptions (A6)-(A7) along with the positiveness of $b''(\cdot)$ from Assumption (A1)  imply that 
for any $\boldsymbol{\beta}_j$ in the interior of $\mathcal{B}$  (and hence in particular for $\boldsymbol{\beta}_j = \boldsymbol{\beta}_j^{M\alpha}$), we have
\begin{eqnarray}
	||\boldsymbol{\beta}_j||_2^2 = O(||\Sigma\boldsymbol{\beta}_0||_2^2) = O(\Lambda_{\max}(\Sigma)) = O(\Lambda_{\max}(\Sigma_1)).
	\label{EQ:eigenCond}
\end{eqnarray}
Here, the first equality is as shown in the proof of Theorem 5 of \cite{Fan/Song:2010} while the remaining equalities are argued above. 

Now, under Assumptions (A1)--(A5), we have the exponential convergence result for the marginal MDPDE as presented in the following Lemma.

\begin{lemma}
	Suppose that (A1)--(A5) hold for a given $\alpha\geq 0$. Then, for any $t>0$,  
	\begin{eqnarray}
		P\left(\sqrt{n}\left|\widehat{\beta}_j^{M\alpha} - \beta_j^{M\alpha}\right|\geq \frac{16k_n^{(\alpha)}}{V}(1+t) \right)
		\leq e^{-\frac{2t^2}{K_n^2}} + nm_1e^{-m_0K_n^{^\tau}},
		~~~~~j=1, \ldots, p,
		\label{EQ:MDPDE_expCov}
	\end{eqnarray}
	where $k_n^{(\alpha)}= (1+\alpha)\left[\frac{m_0}{m_3}K_n^\tau L_\alpha+ |b'(K_nB+B)|L_\alpha + \xi_\alpha(K_nB+B)\right]$.
	\label{LEM:MDPDE_expCov}
\end{lemma}

Note that the constant bounds involved in the above Lemma are independent of the index $j$ leading to 
the uniform convergence of all the marginal regression models through union bound. 
We will utilize this fact to derive the sure screening property of the proposed DPD-SIS 
along with its rate of false positive control (based on (\ref{EQ:eigenCond})),
which is presented in the following theorem.

\begin{theorem}\label{THM:DPD_SIS_Sample}
	Let Assumptions (A1)--(A5) hold for a given $\alpha\geq 0$ and $\frac{n^{1-2\kappa}}{(k_nK_n)^2} \rightarrow \infty$ 
	as $n\rightarrow\infty$, where $k_n = k_n^{(\alpha)}$ is as defined in Lemma \ref{LEM:MDPDE_expCov}. Then the following results hold.
	\begin{itemize}
		\item[(a)] For any given $c_3>0$, there exists $C>0$ such that 
		\begin{eqnarray}
			P\left(\max\limits_{1\leq j \leq p} |\widehat{\beta}_j^{M\alpha} - \beta_j^{M\alpha}| \geq c_3n^{-\kappa}\right)
			\leq p R_n,
		\end{eqnarray}
		where 
		$
		R_n = \left[ e^{-\frac{n^{1-2\kappa}C}{(k_nK_n)^2}} + nm_1e^{-m_0K_n^\tau}\right].
		$
		\item[(b)] If additionally the assumptions of Theorem \ref{THM:DPD_SIS_Pop2} hold, then taking $\gamma_n = c_4 n^{-\kappa}$ with $c_4\leq c_2/2$, 
		we have 
		$$	P\left(\widehat{\mathcal{M}}(\gamma_n) \supset \mathcal{M}_0 \right)\geq 1 - sR_n.$$
		\item[(c)] If additionally Assumptions (A6)--(A7) hold,  taking $\gamma_n = c_5n^{-2\kappa}$, $c_5>0$, we get
		\begin{eqnarray}
			P\left(|\widehat{\mathcal{M}}(\gamma_n)| \leq O(n^{2\kappa}\Lambda_{\max}(\Sigma)) \right)\geq 1 - pR_n.\nonumber
		\end{eqnarray} 
	\end{itemize}
\end{theorem}


It is important to note that the bound $R_n$ in the above theorem is exactly the same 
(except for the value of $k_n=k_n^{(\alpha)}$) as obtained by \cite{Fan/Song:2010} for usual SIS
and it will be exponentially small for standard GLMs with appropriate choices of $K_n$; { see Appendix \ref{APP:Assumption} for detailed discussions}.
Thus, along with the additional robustness property, 
our proposed DPD-SIS at any $\alpha>0$ also enjoys the same optimal rate of convergence and false discovery control
as well as the similar sure screening property as the usual SIS under slightly modified assumptions. 
This is the most striking benefit of our proposal in the context of robust variable screening under high-dimensionality.
Additionally, the sure screening property of the DPD-SIS, as stated in Theorem \ref{THM:DPD_SIS_Sample}(b),
does not depend on the variance and the correlation structure of the covariates for any choices of $\alpha\geq 0$.
However, higher correlation among covariates may surely increase the false positive selection
which can be seen by the dependence of the size of $\widehat{\mathcal{M}}(\gamma_n)$ selected via the DPD-SIS  on $\Sigma$
or more precisely on $\Lambda_{\max}(\Sigma)$ [Theorem \ref{THM:DPD_SIS_Sample}(c)].
As we have less correlation among covariates and hence smaller values of $\Lambda_{\max}(\Sigma)$,
the number of variables selected via our DPD-SIS reduces, leading to less false positives due to its sure independence property.
As in the usual SIS, we can also achieve model selection consistency for DPD-SIS at any $\alpha\geq 0$, i.e., 
$$
P\left(\widehat{\mathcal{M}}(\gamma_n) = \mathcal{M}_0\right) = 1-o(1),
$$ 
under appropriate assumptions on $\Lambda_{\max}(\Sigma)$ along with proper control of $K_n$.
As a particular (extreme) example, it holds with the choice of $\gamma_n$ as in Theorem \ref{THM:DPD_SIS_Sample}(b) 
if we have $\left|\mbox{Cov}(b'(\boldsymbol{X}^T\boldsymbol{\beta}_0), X_j)\right|= o(n^{-\kappa})$ for all $j\notin\mathcal{M}_0$,
along with the other necessary conditions of the theorem depending on $\alpha\geq 0$.

\section{Robust Conditional Screening: The DPD-CSIS}
\label{SEC:DPD-CSIS}

Let us now extend the DPD-SIS approach to conditional screening problems in GLMs.
Suppose that, along with the set-up and notation of Section \ref{SEC:DPD_SIS}, 
information is available to always include a set of $q$ covariates, say $\boldsymbol{X}_{\mathcal{C}}$ { (with $q<n-1$ columns),}
and we need to robustly select variables from the remaining pool of $d:=p-q$ variables (say, $\boldsymbol{X}_{\mathcal{D}}$). 
For simplicity, in this section, we assume no intercept terms, since that can be easily incorporated within the given
$\boldsymbol{X}_{\mathcal{C}}$. Further, without loss of generality, we assume  that $\boldsymbol{X}_{\mathcal{C}}=(X_1, \ldots, X_q)^T$
so that $\boldsymbol{X}_{\mathcal{D}}=(X_{q+1}, \ldots, X_p)^T$; denote $\mathcal{C}=\{1, \ldots, q\}$ and $\mathcal{D}=\{q+1, \ldots, p\}$ and hence
$\boldsymbol{\beta}_\mathcal{C}=(\beta_1, \ldots, \beta_q)^T\in \mathbb{R}^q$
and $\boldsymbol{\beta}_\mathcal{D}=(\beta_{q+1}, \ldots, \beta_p)^T \in \mathbb{R}^d$. 
Now, for a given $\alpha\geq 0$, we may choose the variables from $\boldsymbol{X}_{\mathcal{D}}$ based on the marginal MDPDEs defined as 
\begin{eqnarray}
	\widehat{\boldsymbol{\beta}}_{\mathcal{C}j}^{M\alpha} = \left(\widehat{\boldsymbol{\beta}}_{\mathcal{C}j1}^{M\alpha}, \widehat{\beta}_j^{M\alpha}\right)
	=\arg\min\limits_{\boldsymbol{\beta}_\mathcal{C}, \beta_j} ~~
	\frac{1}{n}\sum_{i=1}^n l_\alpha\left(y_i, \boldsymbol{x}_{i\mathcal{C}}^T\boldsymbol{\beta}_{\mathcal{C}} + \beta_j x_{ij}\right),
	~~~~~j= q+1, \ldots, p,
	\label{EQ:MMDPDE_Cond}
\end{eqnarray}
where $l_\alpha(y, \theta)$ is as defined in Section \ref{SEC:DPD_SIS}, 
and $\boldsymbol{x}_{i\mathcal{C}}$ is the $i$-th observation on $\boldsymbol{X}_\mathcal{C}$.
Then, as in (\ref{EQ:DPD_SIS_set}), given a suitable pre-defined threshold $\gamma_n$, 
we may select the variables in the set 
$\widehat{M}_\alpha(\gamma_n|\mathcal{D}) = \left\{ q+1\leq j \leq p : \left|\widehat{\beta}_j^{M\alpha}\right| \geq \gamma_n  \right\}$.
We refer to this extension as the conditional DPD-SIS,  or the DPD-CSIS in short.
Clearly, the  DPD-CSIS again coincides with the usual CSIS of \cite{Barut/etc:2016} at $\alpha=0$ 
and provides a robust generalization at $\alpha>0$. Further, when the conditioning variable set $\boldsymbol{X}_\mathcal{C}$ is empty 
(or contains only the intercept), we are back to our DPD-SIS. 
We here study the properties of the DPD-CSIS in line with the results derived in Section \ref{SEC:prop_DPDSIS}.
{ 
Note that, throughout this section concerning CSIS, $\mathcal{M}_0$ corresponds to the covariates from $\boldsymbol{X}_\mathcal{D}$ 
having non-zero regression coefficients, and accordingly we now have $ s = |\mathcal{M}_0| < n$. 
}

\subsection{Population-level results: Justifications of DPD-CSIS}

Let us continue with the notation of Section \ref{SEC:prop_DPDSIS_pop} and additionally assume that
$E(X_j|\boldsymbol{X}_\mathcal{C})=0$ for all  $j\in\mathcal{D}$. 
We define the population version (functional) of $\widehat{\boldsymbol{\beta}}_{\mathcal{C}j}^{M\alpha}$ from  (\ref{EQ:MMDPDE_Cond}) as
\begin{eqnarray}
	\boldsymbol{\beta}_{\mathcal{C}j}^{M\alpha} = (\boldsymbol{\beta}_{\mathcal{C}j1}^{M\alpha}, \beta_j^{M\alpha})
	=\arg\min\limits_{\boldsymbol{\beta}_\mathcal{C}, \beta_j} 
	E \left[l_\alpha\left(Y, \boldsymbol{X}_{\mathcal{C}}^T\boldsymbol{\beta}_{\mathcal{C}} + \beta_j X_{j}\right)\right],
	~~~~~j= q+1, \ldots, p.
	\label{EQ:MMDPDEpop_Cond}
\end{eqnarray}
Additionally, let us define the functional for the baseline parameter given only $\boldsymbol{X}_{\mathcal{C}}$, without any additional variable,
as
$\boldsymbol{\beta}_{\mathcal{C}}^{M\alpha} 
=\arg\min\limits_{\boldsymbol{\beta}_\mathcal{C}} 
E \left[l_\alpha\left(Y, \boldsymbol{X}_{\mathcal{C}}^T\boldsymbol{\beta}_{\mathcal{C}} \right)\right]$.
Then, throughout all theoretical discussions of DPD-CSIS, as in \cite{Barut/etc:2016}, 
we need to assume that the functionals $\boldsymbol{\beta}_{\mathcal{C}j}^{M\alpha}$ and $\boldsymbol{\beta}_{\mathcal{C}}^{M\alpha}$ are unique,
i.e., the associated marginal problems are fully identifiable. 
Now, for DPD-CSIS at any given $\alpha\geq 0$,  we consider the random variables $m_{\alpha, j}$, for each $j=q+1, \ldots, p$, defined as
\begin{eqnarray}
	m_{\alpha, j} = \frac{B_\alpha(\boldsymbol{X}_{\mathcal{C}j}^T\boldsymbol{\beta}_{\mathcal{C}j}^{M\alpha}) - 
		B_\alpha(\boldsymbol{X}_{\mathcal{C}}^T\boldsymbol{\beta}_{\mathcal{C}}^{M\alpha})}{
		\boldsymbol{X}_{\mathcal{C}j}^T\boldsymbol{\beta}_{\mathcal{C}j}^{M\alpha}- \boldsymbol{X}_{\mathcal{C}}^T\boldsymbol{\beta}_{\mathcal{C}}^{M\alpha}},
	\label{EQ:m_alphaj}
\end{eqnarray} 
where $\boldsymbol{X}_{\mathcal{C}j}=(\boldsymbol{X}_{\mathcal{C}}^T, X_j)^T$ for each $j$ and $B_\alpha$ is as defined in Section \ref{SEC:prop_DPDSIS_pop}.
Denote by $\mathcal{M}_{0\mathcal{D}} = \mathcal{M}_0 \cap \mathcal{D}$ the indices of the truly important variables in $\boldsymbol{X}_\mathcal{D}$.
Then, we have the following results, in analogue of Theorems \ref{THM:DPD_SIS_Pop1} and \ref{THM:DPD_SIS_Pop2}, 
that justify our DPD-CSIS algorithm as a reasonable procedure for conditional screening. 
Here, in analogue of \cite{Barut/etc:2016}, we define 
$\mbox{Cov}_L(Y, X_j|\boldsymbol{X}_{\mathcal{C}})
:= E\left[\left(Y-E_L[Y|\boldsymbol{X}_{\mathcal{C}}]\right)\left(X_j-E_L[X_j|\boldsymbol{X}_{\mathcal{C}}]\right)\right]$,
for any $j\in\mathcal{D}$, where $E_L[\cdot|\boldsymbol{X}_\mathcal{C}]$ denote the best linear regression fit given $\boldsymbol{X}_\mathcal{C}$;
clearly $E_L[Y|\boldsymbol{X}_{\mathcal{C}}]=b'(\boldsymbol{X}_{\mathcal{C}}^T\boldsymbol{\beta}_{\mathcal{C}}^{M\alpha}).$

\begin{theorem}\label{THM:DPD_SIS_Pop1_Cond}
	For a given $\alpha\geq 0$ and  any $j\in\mathcal{D}$, the (conditional) 
	marginal  MDPDE functional $\beta_j^{M\alpha}$ in  (\ref{EQ:MMDPDEpop_Cond}) is zero if and only if 
	$\mbox{Cov}_L(Y, X_j|\boldsymbol{X}_{\mathcal{C}})=0$.
\end{theorem}

\begin{theorem}\label{THM:DPD_SIS_Pop2_Cond}
	Given any $\alpha\geq 0$, suppose that $E[m_{\alpha,j}X_j^2]\leq c_2$ uniformly in $j\in \mathcal{D}$,  for some constant $c_2$.
	If there exist constants $c_1>0, \kappa<-1/2$ such that  
	$\left|\mbox{Cov}_L(Y, X_j|\boldsymbol{X}_\mathcal{C})\right|\geq c_1n^{-\kappa}$ for all $ j\in\mathcal{M}_{0\mathcal{D}}$, 
	then we have 
	$
	\min\limits_{j\in\mathcal{M}_{0\mathcal{D}}}|\beta_j^{M\alpha}| \geq c_3 n^{-\kappa}, 
	$ 
	for another constant $c_3>0$.
\end{theorem}


\subsection{Sample-level properties: Sure Screening via DPD-CSIS}

We now extend the results of Section \ref{SEC:prop_DPDSIS_sample} for the unconditional DPD-SIS to the case of conditional screening to show the uniform convergence of the associated (conditional) MDPDEs and the resulting sure screening property 
of the DPD-CSIS. We continue with the notation of Section  \ref{SEC:prop_DPDSIS_sample}
and assume that Assumptions (A1)--(A7) hold with $\boldsymbol{\beta}_j$ and $\boldsymbol{\beta}_j^{M\alpha}$
replaced by $\boldsymbol{\beta}_{\mathcal{C}j}\in\mathbb{R}^{q+1}$ and $\boldsymbol{\beta}_{\mathcal{C}j}^{M\alpha}$, respectively, in (A2)--(A4).
We also assume the following additional condition.

\begin{itemize}
	\item[(A8)] There exists  $C>0$ such that $\Lambda_{\min}\left(E\left[m_{\alpha,j} \boldsymbol{X}_{\mathcal{C}j}\boldsymbol{X}_{\mathcal{C}j}^T\right]\right)>C$, uniformly over $j\in\mathcal{D}$. 
	
\end{itemize} 

Note that Assumption (A8) is mild (and regular) if $B_\alpha$ is strictly convex implying $m_{\alpha,j}>0$ almost surely.
Further, we define
$\boldsymbol{Z} = E\left(E_L[\boldsymbol{X}_{\mathcal{D}}|\boldsymbol{X}_{\mathcal{C}}]\left[\boldsymbol{X}^T\boldsymbol{\beta}_0
- \boldsymbol{X}_{\mathcal{C}}^T\boldsymbol{\beta}_{\mathcal{C}}^{M\alpha}\right]\right)$ and 
$\Sigma_{\mathcal{D}|\mathcal{C}} = E\left(\boldsymbol{X}_{\mathcal{D}} - E_L[\boldsymbol{X}_{\mathcal{D}}|\boldsymbol{X}_{\mathcal{C}}]\right)
\left(\boldsymbol{X}_{\mathcal{D}} - E_L[\boldsymbol{X}_{\mathcal{D}}|\boldsymbol{X}_{\mathcal{C}}]\right)^T$.
We can show  that Assumptions (A6)-(A8) imply the following analogue of (\ref{EQ:eigenCond}) for this conditional case, given by
\begin{eqnarray}
	||\boldsymbol{\beta}_{\mathcal{D}}||_2^2 = O\left(\Sigma_{\mathcal{D}|\mathcal{C}}+\boldsymbol{Z}\boldsymbol{Z}^T\right).
	\label{EQ:eigenCond0}
\end{eqnarray}
Then, we have the desired results in analogue to Theorem \ref{THM:DPD_SIS_Sample} for the present conditional case of the DPD-CSIS
which is presented in the following theorem. The proof is similar to that of Theorem \ref{THM:DPD_SIS_Sample}, 
but using (\ref{EQ:eigenCond0}) instead of (\ref{EQ:eigenCond}), and is hence omitted for brevity. 

\begin{theorem}\label{THM:DPD_SIS_Sample-Cond}
	Suppose that, for a given $\alpha\geq 0$, Assumptions (A1)--(A5) hold with $\boldsymbol{\beta}_j$ and $\boldsymbol{\beta}_j^{M\alpha}$
	replaced by $\boldsymbol{\beta}_{\mathcal{C}j}\in\mathbb{R}^{q+1}$ and $\boldsymbol{\beta}_{\mathcal{C}j}^{M\alpha}$, respectively, in (A2)--(A4). 
	Also, let $\frac{n^{1-2\kappa}}{(k_nK_n)^2} \rightarrow \infty$ as $n\rightarrow\infty$, 
	where $k_n=k_n^{(\alpha)}$ is as defined in Lemma \ref{LEM:MDPDE_expCov}. Then, the following results hold.
	\begin{itemize}
		\item[(a)] For any given $c_3>0$, there exists $C>0$ such that 
		\begin{eqnarray}
			P\left(\max\limits_{q+1\leq j \leq p} |\widehat{\beta}_j^{M\alpha} - \beta_j^{M\alpha}| \geq c_3n^{-\kappa}\right)
			\leq d R_n,
		\end{eqnarray}
		where $R_n$ is as defined in Theorem \ref{THM:DPD_SIS_Sample}.
		\item[(b)] If additionally the assumptions of Theorem \ref{THM:DPD_SIS_Pop2_Cond} hold, 
		then taking $\gamma_n = c_4 n^{-\kappa}$ with $c_4\leq c_2/2$, 
		we have 
		$$		P\left(\widehat{\mathcal{M}}(\gamma_n) \supset \mathcal{M}_0 \right)\geq 1 - sR_n.$$
		\item[(c)] If additionally Assumptions (A6)--(A8) hold, taking $\gamma_n = c_4n^{-2\kappa}$, $c_4>0$, we get
		\begin{eqnarray}
			P\left(|\widehat{\mathcal{M}}(\gamma_n)| \leq 
			O\left(n^{2\kappa}\Lambda_{\max}\left(\Sigma_{\mathcal{D}|\mathcal{C}}+\boldsymbol{Z}\boldsymbol{Z}^T\right)\right)
			\right) \geq 1 - dR_n.
			\label{EQ:FDR_Cond}
		\end{eqnarray} 
	\end{itemize}
\end{theorem}

Note that the rate of convergence in the above theorem is exactly the same as in the unconditional case 
(Theorem \ref{THM:DPD_SIS_Sample}) and that they are in line with the existing literature on variable screening. 
In the particular case of the linear regression model, we have $\boldsymbol{Z}=\boldsymbol{0}$, 
and hence the result (\ref{EQ:FDR_Cond}) in Theorem \ref{THM:DPD_SIS_Sample-Cond} reduces to 
\begin{eqnarray}
	P\left(|\widehat{\mathcal{M}}(\gamma_n)| \leq 
	O\left(n^{2\kappa}\Lambda_{\max}\left(\Sigma_{\mathcal{D}|\mathcal{C}}\right)\right)
	\right) \geq 1 - dR_n.
	\label{EQ:FDR_Cond0}
\end{eqnarray} 
In general, if we additionally assume $||\boldsymbol{Z}||_2^2 = o\left(\Lambda_{\min}\left(\Sigma_{\mathcal{D}|\mathcal{C}}\right)\right)$, 
as in Condition 3(iii) of \cite{Barut/etc:2016}, we can also have (\ref{EQ:FDR_Cond0}) instead of (\ref{EQ:FDR_Cond}) in 
Theorem \ref{THM:DPD_SIS_Sample-Cond}.

{
\section{Robustness Property: Theoretical Justifications}
\label{SEC:IF}

The robustness of the proposed DPD-SIS and DPD-CSIS under data contamination follows directly from the robustness of the associated marginal MDPDEs $\widehat{\boldsymbol{\beta}}_j^{M\alpha}$.
For an intuitive understanding, recall that the MDPDE estimating equation downweights the outliers with a weight $f(y; \theta)^{\alpha}$ to achieve robustness.
As the value of $\alpha>0$ increases, more down-weighting takes place, reducing the contribution of outliers in the estimation process, 
which leads to improved robustness of the MDPDEs and the subsequent DPD-SIS or DPD-CSIS procedure. 
At $\alpha=0$, there is no down-weighting of the outlying observations and so the resulting MLE or the associated SIS/CSIS are non-robust under data contamination.

The robustness characteristic of both DPD-SIS and DPD-CSIS, i.e., their increasing robustness with increasing $\alpha>0$, 
compared to the usual SIS or CSIS (at $\alpha=0$) can be theoretically justified by the classical influence function analyses 
under the Huber's $\epsilon$-contamination model \citep{Hampel/etc:1986}.
The influence function (IF) provides a measure of asymptotic bias, in any statistical functional, caused by infinitesimal contamination at a distant outlying point.
If this IF tends to infinity (with either sign) as the contamination point moves further away, the resulting bias then increases indefinitely under contamination
indicating the non-robust nature of the associated functional (e.g., estimator).
However, as long as the IF remains bounded as a function of the contamination point, 
the resulting functional (estimator) cannot have a value extremely far from the true value even under (infinitesimal) contamination at a very distant point,
which justifies its robustness; the smaller the maximum extent of the IF (in absolute value) the greater the stability of the associated estimator.   
Therefore, the IF of the marginal MDPDEs would then give a theoretical justification of their robustness, 
and hence, the same for the proposed DPD-SIS and DPD-CSIS as well.
 
The existing theory of the MDPDE \citep{Basu/etc:2011,Ghosh/Basu:2013,Ghosh/Basu:2016} has covered its IF under different parametric models including the GLMs. 
In particular, the IF of the marginal MDPDE $\widehat{\boldsymbol{\beta}}_j^{M\alpha}$ in the present context of DPD-SIS under a given GLM
would have the form \citep{Ghosh/Basu:2016}
\begin{eqnarray}
	\mathcal{IF}_j^{(\alpha)}(y_t|x_{jt}) &=& \boldsymbol{\Psi}_n^{-1} \cdot \psi_\alpha(y_t, \beta_{j0} + \beta_j x_{jt})[1, x_{jt}]^T,
	\nonumber
\end{eqnarray}
where $y_t$ is the contamination point in the response variable with the associated covariate value being $x_{jt}$,
$\beta_{j0}$ and  $\beta_j$ are the assumed true parameter values, $\psi_\alpha$ is as defined in (\ref{EQ:DPD_psi}), 
and $\boldsymbol{\Psi}_n$ is some suitable matrix independent of the contamination point 
(see Eq.~(5) of \cite{Ghosh/Basu:2016} for its exact form). 
For DPD-CSIS, the IF of the corresponding marginal MDPDE, defined in (\ref{EQ:MMDPDE_Cond}), has the form 
\begin{eqnarray}
	\mathcal{IF}_j^{(\alpha)}(y_t|x_{jt}) &=& \boldsymbol{\Psi}_n^{-1} \cdot \psi_\alpha(y_t, \boldsymbol{x}_{t\mathcal{C}}^T\boldsymbol{\beta}_{\mathcal{C}} 
	+ \beta_j x_{jt})[\boldsymbol{x}_{t\mathcal{C}}^T, x_{jt}]^T,
	\nonumber
\end{eqnarray}
where $\boldsymbol{x}_{t\mathcal{C}}$ denotes the values of the conditioning covariates associated with the contaminated response $y_t$. 
Therefore, in both cases, we can see that the form of function $\psi_\alpha$ determines the boundedness of the IFs of the marginal MDPDEs,
and hence, the robustness of the resulting DPD-SIS and DPD-CSIS procedure. 
But, it can be easily noted that the function $\psi_\alpha$, given in (\ref{EQ:DPD_psi}), is bounded for all standard GLMs at any $\alpha>0$  
and is unbounded at $\alpha=0$ for most GLMs (specifically where the response has an unbounded support).
Thus, the proposed DPD-SIS and DPD-CSIS with any $\alpha>0$ would be robust under all GLMs.
Further, the maximum value of $|\psi_\alpha|$  also decreases as $\alpha$ increases, indicating the increasing robustness of the associated marginal MDPDEs
for increasing values of $\alpha>0$; the same is also transferred subsequently for the proposed DPD-SIS and DPD-CSIS procedures with any $\alpha>0$.
Note that, this analysis additionally yields a theoretical justification of the non-robustness of the usual SIS at $\alpha=0$ (with unbounded IF)
for all GLMs having unbounded support for the response distribution.

The above-mentioned IFs can be further investigated for specific GLMs, by looking at the corresponding $\psi_\alpha$ function.
For the linear regression models, for example, the associated $\psi_\alpha(y_t, \theta)$ has a form proportional to 
$\left(y_t - \theta\right)\exp\left({-\frac{\alpha\left(y_t - \theta\right)^2}{2(\sigma_j^{(0)})^2}}\right)$;
this case has been studied extensively in \cite[][Appendix A.1]{Ghosh/Thoresen:2020}.
For the case of Poisson regression models, the $\psi_\alpha$ function has the form 
\begin{eqnarray}
\psi_\alpha(y_t, \theta)=  \frac{(y_t - e^\theta)}{(y_t!)^\alpha}\exp(\alpha y_t \theta + \alpha e^\theta) - \xi_\alpha(\theta).
\nonumber
\end{eqnarray}
Like the case of linear regression, one can also here rigorously see that, for any given covariate values (and hence any given $\theta$)
the function $\psi_\alpha(y_t, \theta)$ is bounded in $y_t$ for all $\alpha>0$, unbounded at $\alpha=0$ 
and $\sup_{y_t} |\psi_\alpha(y_t, \theta)|$ decreases as $\alpha>0$ increases. 
So, all the general robustness properties of the marginal MDPDEs, and  hence, those for the DPD-SIS and DPD-CSIS,  clearly hold also for Poisson regression models. 
For another important GLM, logistic regression, the $\psi_\alpha$ function, and hence, 
the IF of the MDPDEs would be bounded for all $\alpha\geq 0$ due to the bounded support of $y_t$;
but its maximum value would still decrease indicating greater extent of robustness for DPD-SIS or DPD-CSIS with increasing values of $\alpha$.


}

{

\section{Numerical Illustrations}

To illustrate the finite-sample performance of the proposed DPD-SIS and DPD-CSIS, 
we have performed extensive simulation studies for several important examples of GLMs;
for brevity, a few interesting cases of the linear and logistic regression models are presented in this section.
The corresponding \texttt{R} codes (for linear, logistic and also Poisson regressions) are provided in a public GitHub repository 
titled \texttt{dpdSIS}  (available at \url{https://github.com/abhianik/dpdSIS}). 

We would like to mention that prior numerical illustrations on the performance of the DPD-SIS under linear regression models, 
along with its comparison with several other existing robust SIS procedures, are available in \cite{Ghosh/Thoresen:2020}.
However, there is no literature on robust variable screening procedures for general GLMs (beyond linear regression), 
except for the rank-correlation based SIS (rank-SIS) of \cite{Li/etc:2012a}.
So,  in the present paper,  we have compared the DPD-SIS and DPD-CSIS, under the logistic regression model,  
with the classical SIS of \cite{Fan/Song:2010} and the only existing robust rank-SIS of \cite{Li/etc:2012a}, 
although rank-SIS is a non-parametric method and thus philosophically different from our approach.
It is also worth pointing out that even though \cite{Li/etc:2012a} did use their rank-SIS method for logistic models, 
the properties of their approach is proven only for linear models unlike our proposed DPD-SIS and DPD-CSIS whose sure screening properties are proven theoretically for general GLMs including the logistic regression.
Finally, for consistency, the illustrations with linear regression models are also made only with the usual SIS/CSIS and rank-SIS in the present paper 
(which is indeed sufficient in view of the existing comparisons in \cite{Ghosh/Thoresen:2020} for linear regression models).

\subsection{Simulation Settings}\label{SEC:SimSettings}
}

We simulate each sample of covariates (except intercept) from a multivariate normal distribution with mean vector $\boldsymbol{0}$ 
and a variance matrix having $(i,j)$-th element as $\rho$ for all $i\neq j$ and 1 for $i=j$; 
clearly $\rho=0$ yields the case of independent covariates whereas a non-zero value of $\rho$ indicates correlated covariates.    
The intercept term (1) is then added as the first covariate. 
{ Then the responses are generated according to a specified GLM (linear or logistic)} with (true) coefficient value 
$\boldsymbol{\beta}_0=(\beta_{01}, \ldots, \beta_{0p})$.
Considering $s=4$ (i.e., four covariates are actually related to the response), 
we choose $\boldsymbol{\beta}_0$ such that four components (in addition to the first one, the intercept) are non-zero and the rest are zero. 
{ In particular, non-zero coefficients of $\boldsymbol{\beta}_0$ are considered to be sparsely distributed in positions 1 (intercept), 2, 6, 26 and 126. 
}
Different values of $\boldsymbol{\beta}_0$ in terms of both the position and strength of these four non-zero coefficients
are considered along with different values of $\rho, n$ and $p$; the results obtained for $p=5000$, $n=100$ and $\rho=0, 0.3$ are presented here. 
{ For linear regression, errors are generated from the standard normal distribution and the error variance is assumed to be known
(in consistence with the theoretical set-up of the present paper).}

For each scenario,  the usual SIS and the proposed DPD-SIS are  applied to the simulated sample of size $n=100$
to identify the top $n-1=99$ variables; it is then examined if the true non-zero (significant) covariates are selected.
The full process is repeated 300 times and the number of true positives 
(number of selected variables having true non-zero coefficients)  are studied as a summary measure. 
{ As mentioned previously, for comparison, we repeat the same exercise for the usual SIS of \cite{Fan/Song:2010} and 
the existing robust rank-SIS of \cite{Li/etc:2012a} and the corresponding results are also presented in the same figures. }

Additionally, to examine the robustness, we repeat all these simulations again by contaminating 10\% of the observations in 150 non-significant covariates 
(having true regression coefficient zero) by independent observations (outliers) generated from a normal distribution with mean equal $-10$ and variance equal 1.
{ In order to have a clear outlier-effect, contaminations are introduced to observations having smaller responses in case of linear regression,
and smaller probability of success for logistic model.
Several other variations of the contamination schemes are also studied 
(e.g., introducing contaminations to randomly chosen observations, or using different contamination distributions, etc.) 
but they all produce similar results leading to the exact same conclusions about the performances of the DPD-SIS in comparison to the usual SIS and rank-SIS;
so, these results are not reported here for brevity.} 

\begin{figure}[!h]
	\centering
	\subfloat[Independent Covariates ($\rho = 0$), Strong Signal (non-zero values of $\boldsymbol{\beta}_{0}$ are all 5)]{
		\includegraphics[width=0.425\textwidth]{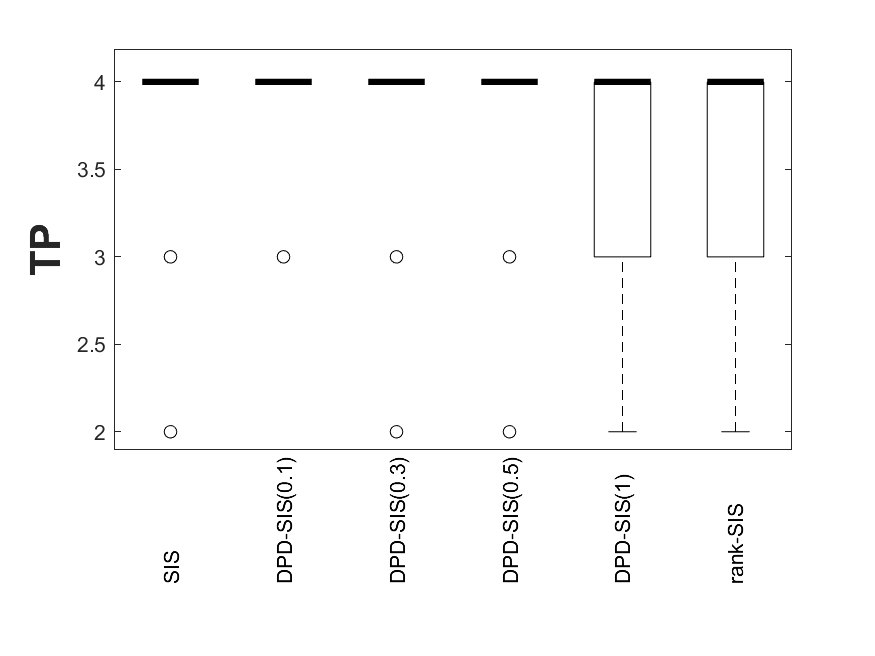}
		\includegraphics[width=0.425\textwidth]{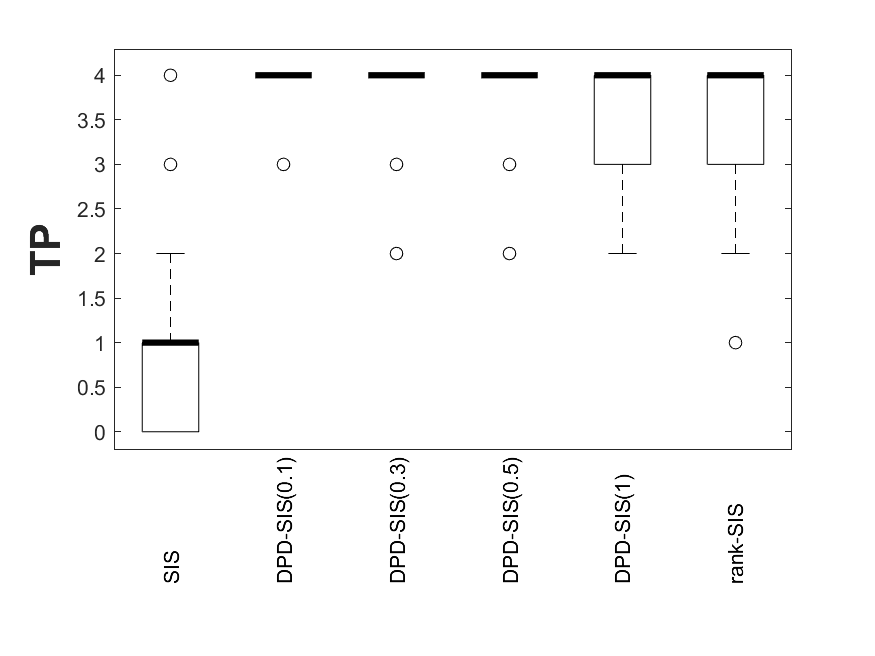}
		\label{FIG:TP_logistic_A}}\\			
	\subfloat[Independent Covariates, Weaker Signal $(\beta_{01}, \beta_{02}, \beta_{06}, \beta_{0,26}, \beta_{0,126})=(1,2,3,1,5)$]{
		\includegraphics[width=0.425\textwidth]{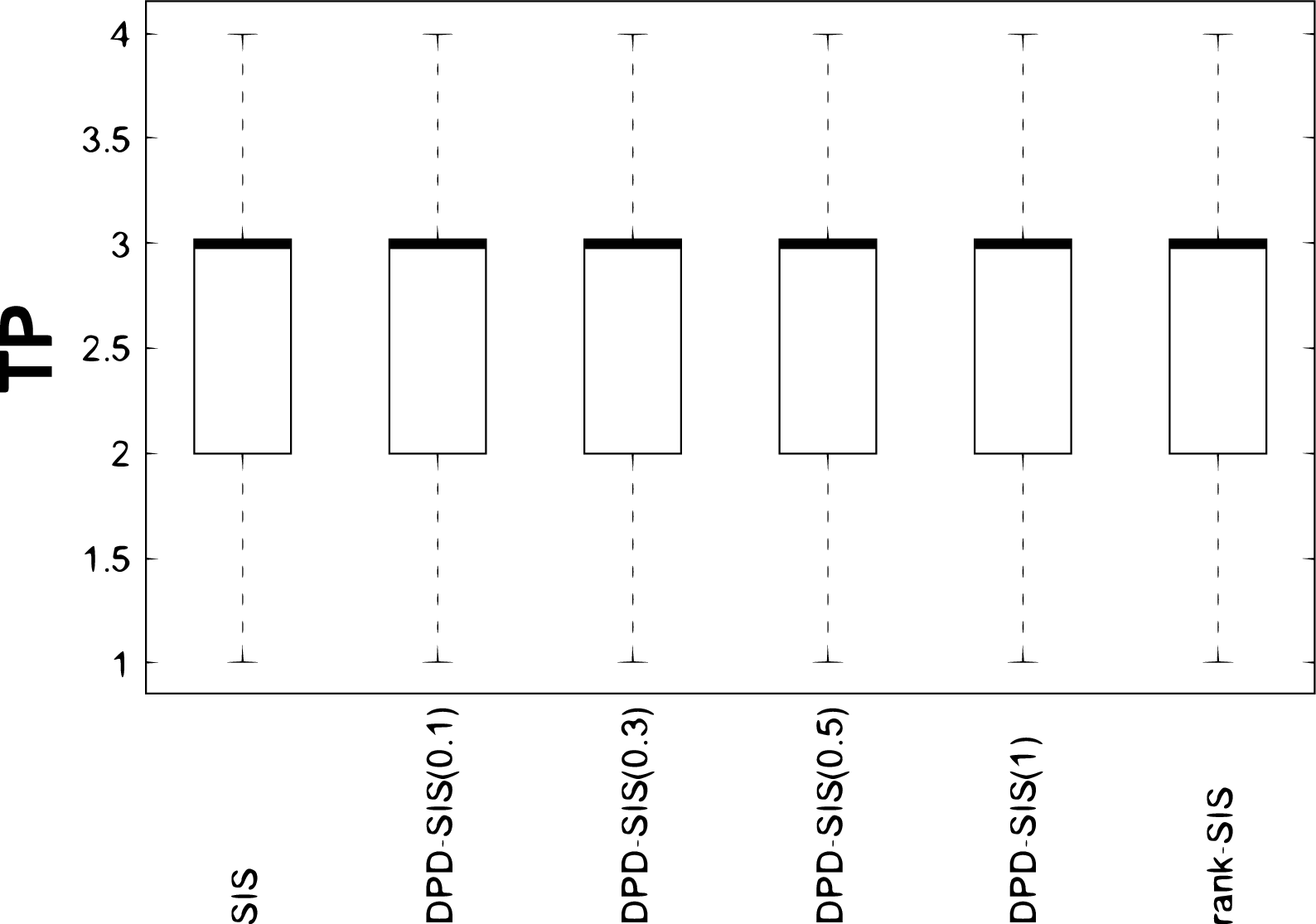}
		\includegraphics[width=0.425\textwidth]{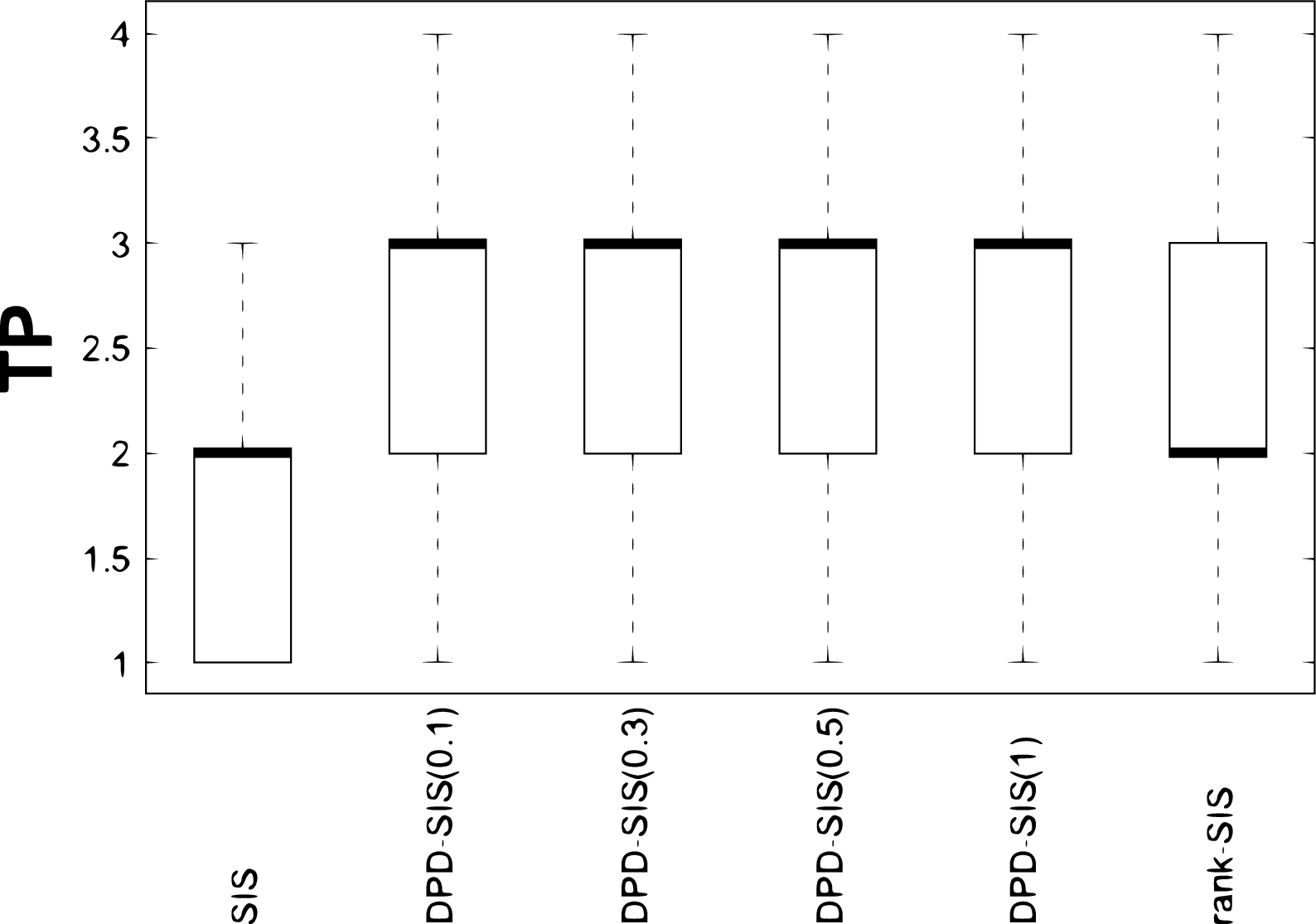}
		\label{FIG:TP_logistic_B}}\\			
	\subfloat[Dependent Covariates ($\rho = 0.3$), Strong Signal  (non-zero values of $\boldsymbol{\beta}_{0}$ are all 5)]{
		\includegraphics[width=0.425\textwidth]{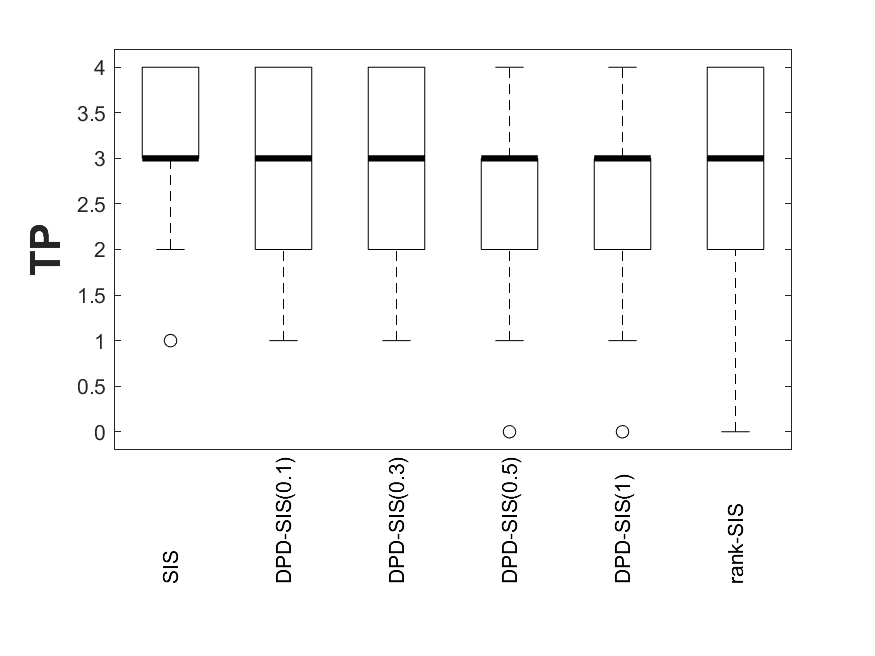}
		\includegraphics[width=0.425\textwidth]{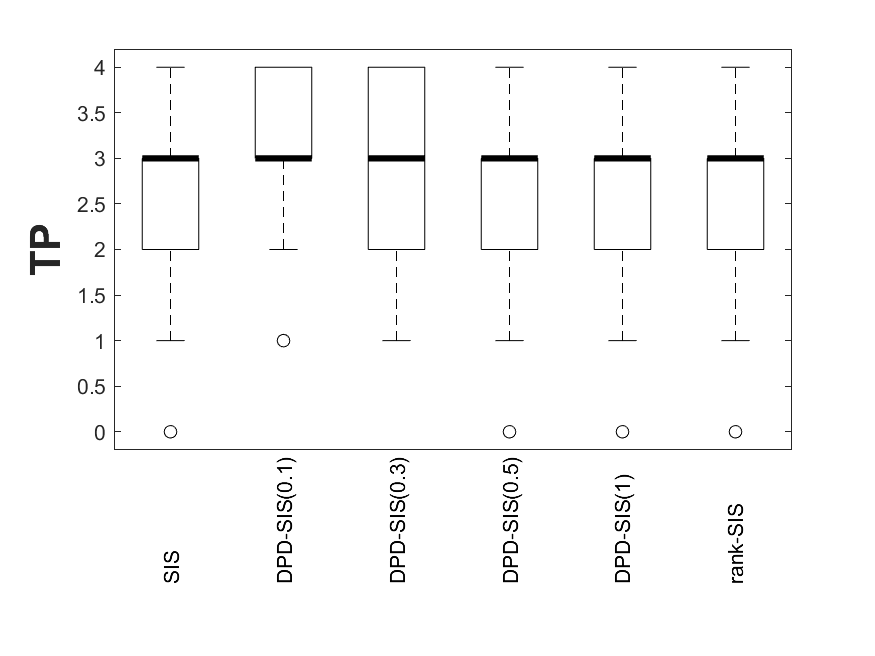}
		\label{FIG:TP_logistic_C}}
	\caption{The box-plots of true-positives selected by the DPD-SIS at different $\alpha$ for different simulation set-ups { under the logistic regression model} 
		with pure data (left panel) and 10\% contaminated data (right panel). 
		{ The results for usual SIS of \cite{Fan/Song:2010} and rank-SIS of \cite{Li/etc:2012a} are also presented in the same plots for comparisons.} }
	\label{FIG:TP_logistic}
\end{figure}

\begin{figure}[!h]
	\centering
	\subfloat[Independent Covariates ($\rho = 0$), Strong Signal (non-zero values of $\boldsymbol{\beta}_{0}$ are all 5)]{
		\includegraphics[width=0.425\textwidth]{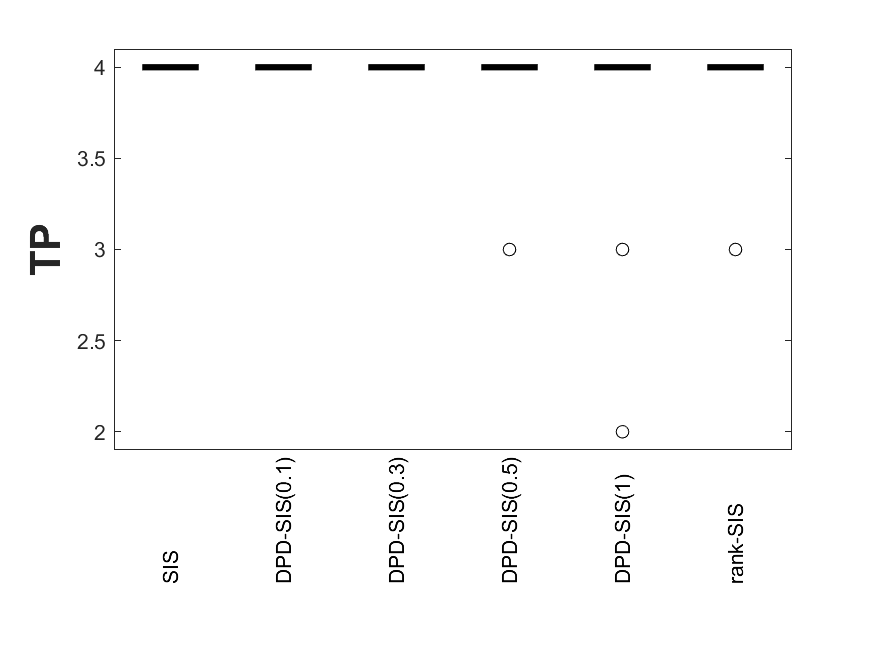}
		\includegraphics[width=0.425\textwidth]{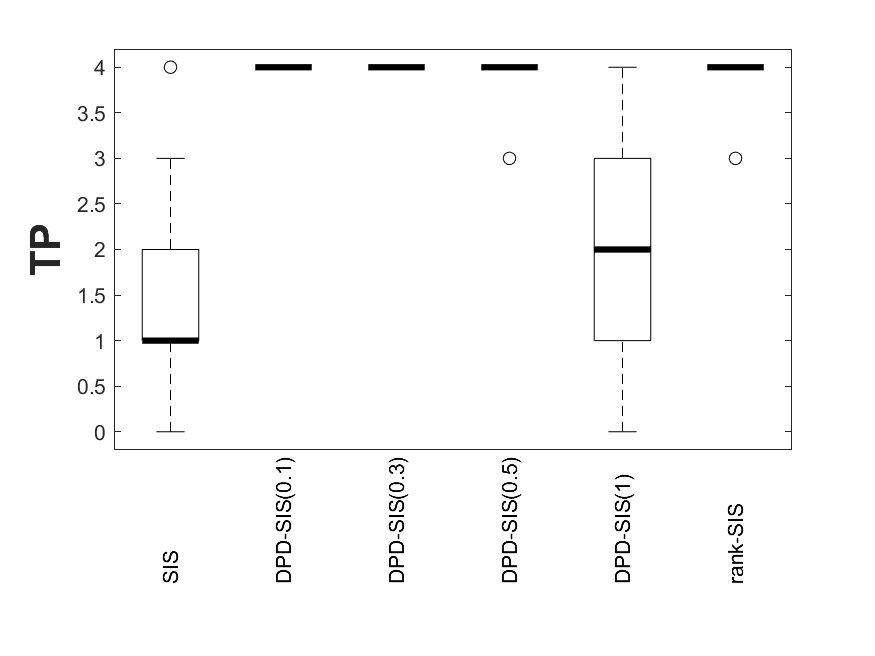}
		\label{FIG:TP_lrm_A}}\\			
	\subfloat[Independent Covariates, Weaker Signal $(\beta_{01}, \beta_{02}, \beta_{06}, \beta_{0,26}, \beta_{0,126})=(1,2,3,1,5)$]{
		\includegraphics[width=0.425\textwidth]{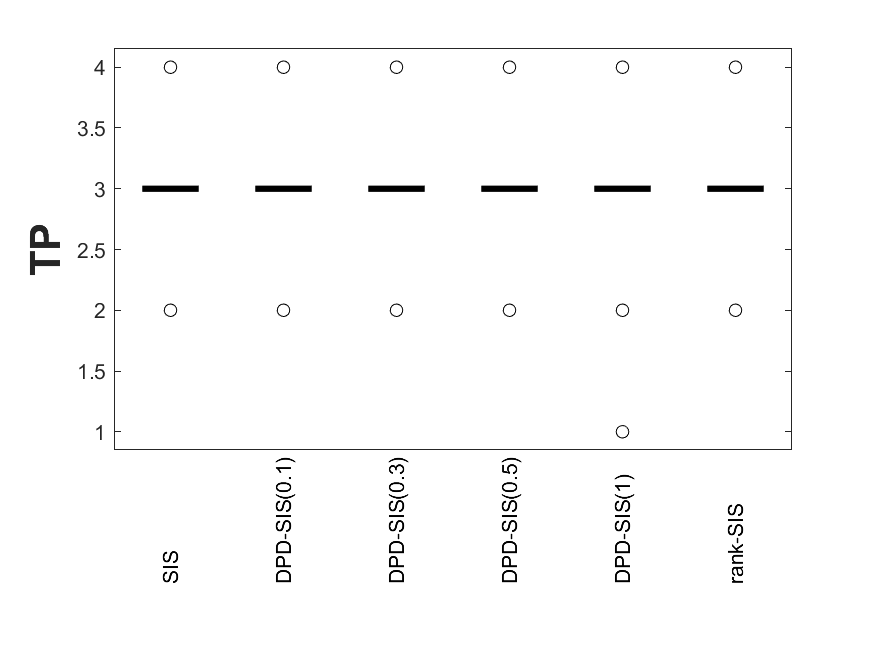}
		\includegraphics[width=0.425\textwidth]{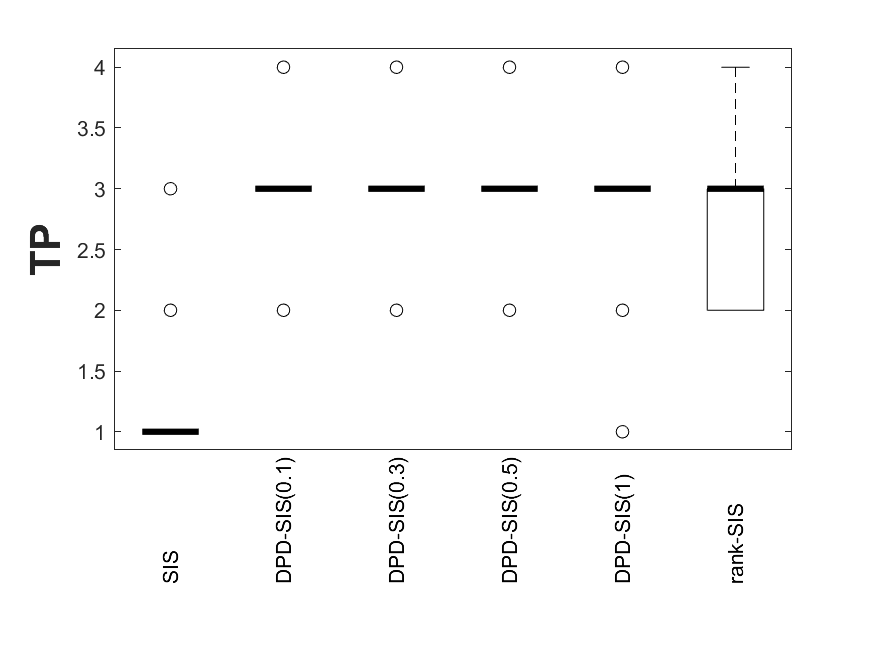}
		\label{FIG:TP_lrm_B}}\\			
	\subfloat[Dependent Covariates ($\rho = 0.3$), Strong Signal  (non-zero values of $\boldsymbol{\beta}_{0}$ are all 5)]{
		\includegraphics[width=0.425\textwidth]{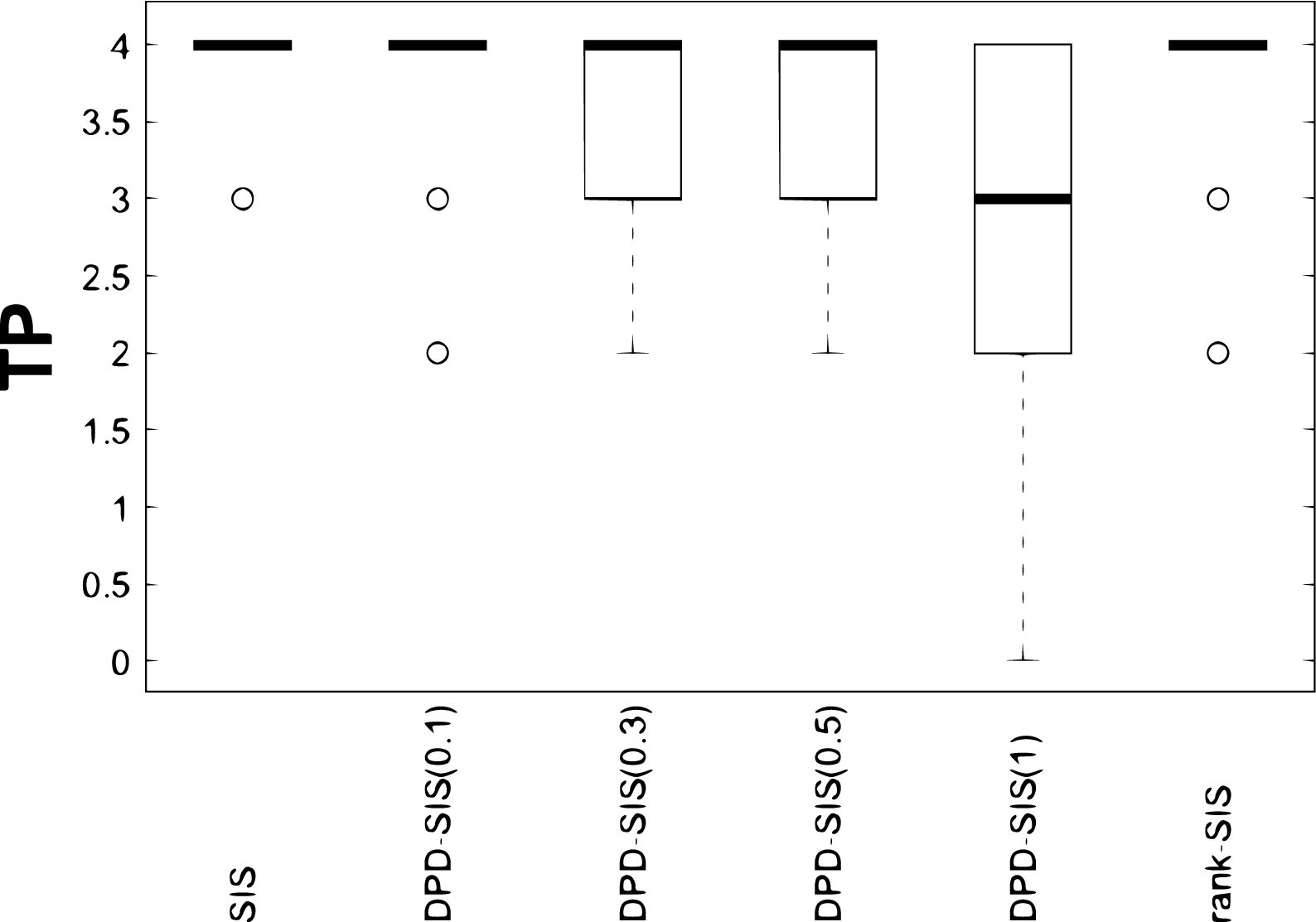}
		\includegraphics[width=0.425\textwidth]{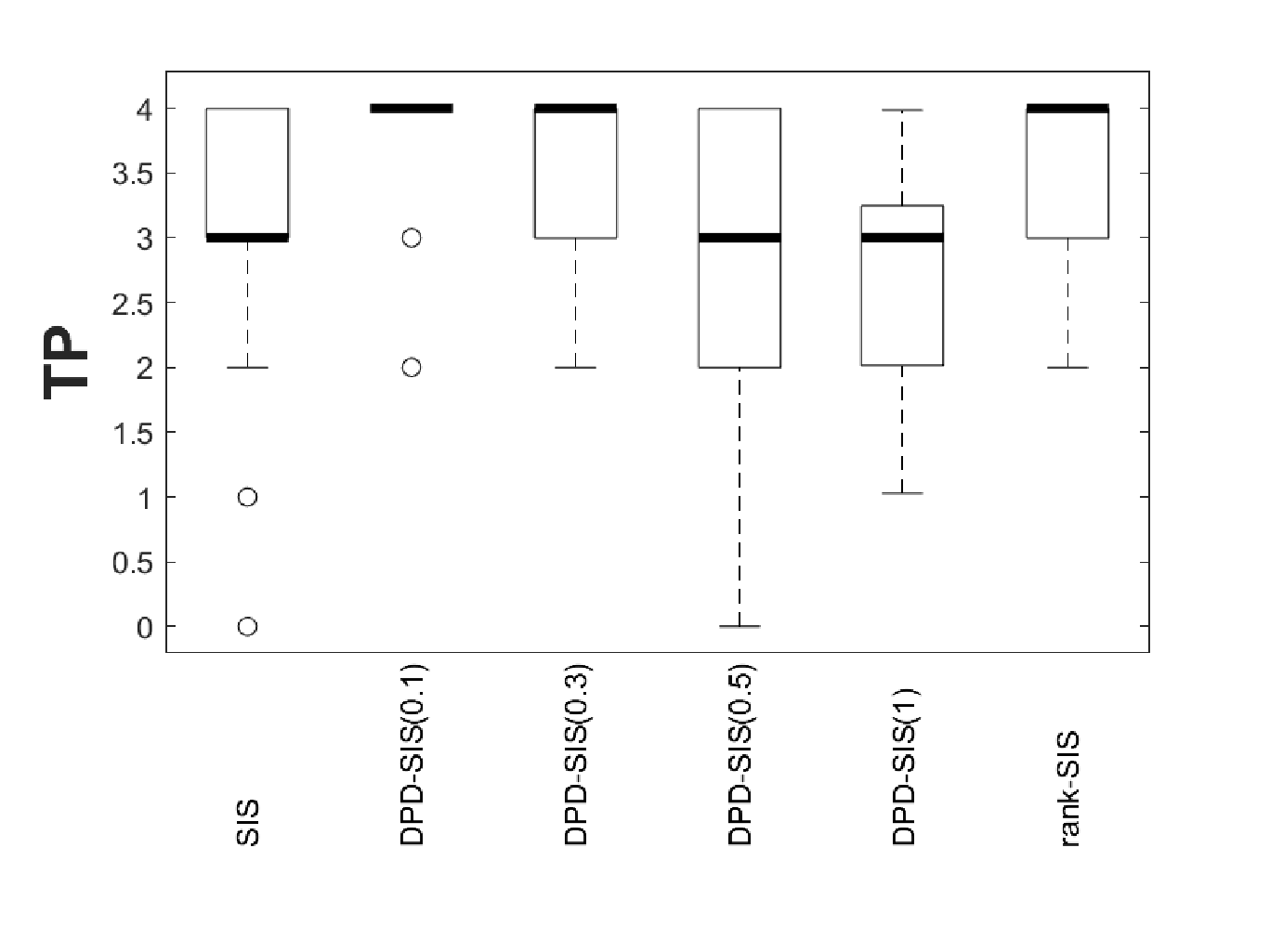}
		\label{FIG:TP_lrm_C}}
	\caption{{ The box-plots of true-positives selected by the DPD-SIS at different $\alpha$ for different simulation set-ups under the linear regression model 
		with pure data (left panel) and 10\% contaminated data (right panel). 
		 The results for usual SIS of \cite{Fan/Song:2010} and rank-SIS of \cite{Li/etc:2012a} are also presented in the same plots for comparisons.} }
	\label{FIG:TP_lrm}
\end{figure}

{

\subsection{Simulation Results: Performance of DPD-SIS}

We present the box-plots of true-positives obtained by the proposed DPD-SIS procedure at different $\alpha>0$ 
(along with the usual SIS and the rank-SIS) both without and with data contamination (as specified in Section \ref{SEC:SimSettings})
in Figures \ref{FIG:TP_logistic} and \ref{FIG:TP_lrm} for the logistic and the linear regression cases, respectively.
}

We observe that for independent covariates and relatively strong signal (larger coefficient values)
the DPD-SIS with smaller $\alpha\leq 0.5$ performs exactly similar to the usual SIS under pure data. 
But, when there is contamination in the data, the performance of the usual SIS deteriorates significantly 
whereas the proposed DPD-SIS remains stable and yields better variable selection results, ignoring the effect of outliers
(see, e.g., Figure \ref{FIG:TP_logistic_A}). 
As the signal gets weaker, the DPD-SIS fails to select all true positives but it still performs as good as the usual SIS under pure data;
in case of additional contamination, the DPD-SIS remains much more stable even with weaker signal, 
although the usual SIS gets significantly affected (see, e.g., Figure \ref{FIG:TP_logistic_B}).  
However, the performance of the proposed DPD-SIS as well as the usual SIS becomes significantly worse when the covariates are strongly correlated
which is expected from our theory as well. In such a case, the effect of outliers may not be so prominent as in the case of independent covariates,
but there is still a slight decrease in the number of true positives obtained by the usual SIS under contamination
and DPD-SIS again performs robustly as claimed (see, e.g., Figure \ref{FIG:TP_logistic_C}).

{ 
The rank-SIS method performs better than the usual SIS in terms of robustness under data contamination but worse under pure data.
Additionally, the results obtained by the rank-SIS are clearly worse than the best results obtained by our proposed DPD-SIS at a suitable $\alpha$  under both pure and contaminated data scenarios,
justifying the advantages of our proposal over the existing rank-SIS method under both linear and logistic regression models.  
}

\begin{figure}[!h]
	\centering
	\subfloat[Logistics Regression Models]{
		\includegraphics[width=0.425\textwidth]{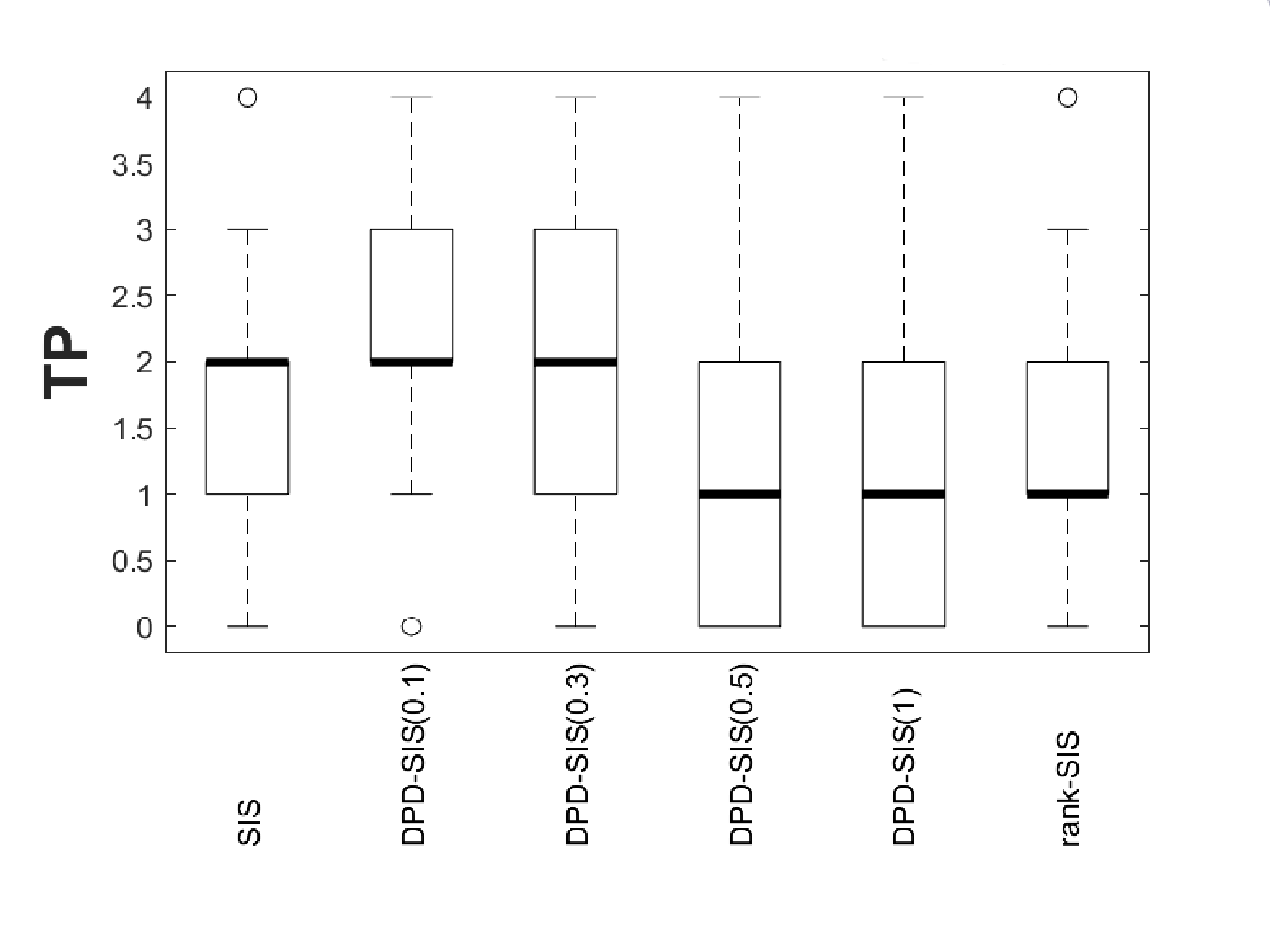}
		\includegraphics[width=0.425\textwidth]{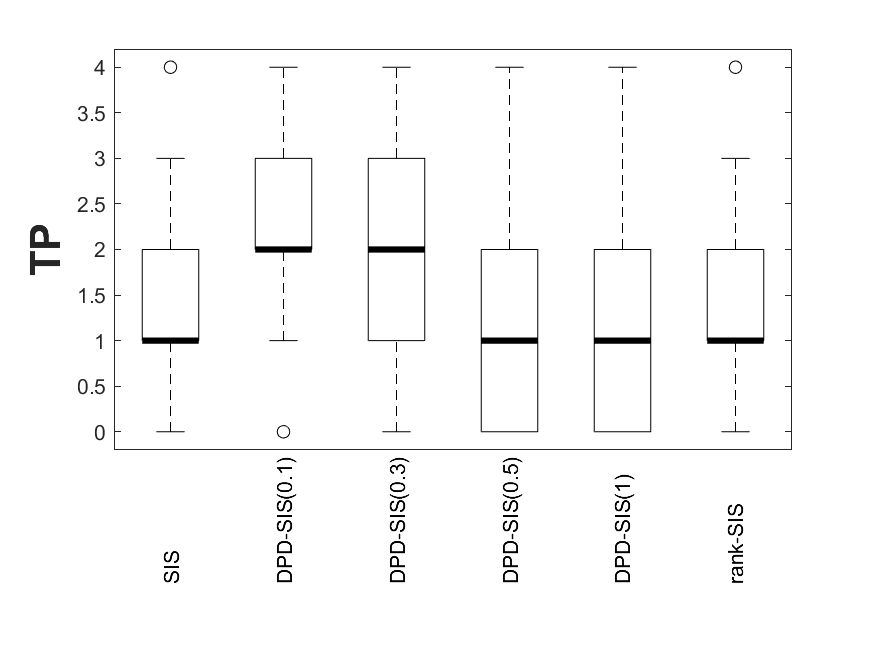}
		\label{FIG:TP_logistic_CSIS}}\\			
	\subfloat[Linear Regression Models]{
		\includegraphics[width=0.425\textwidth]{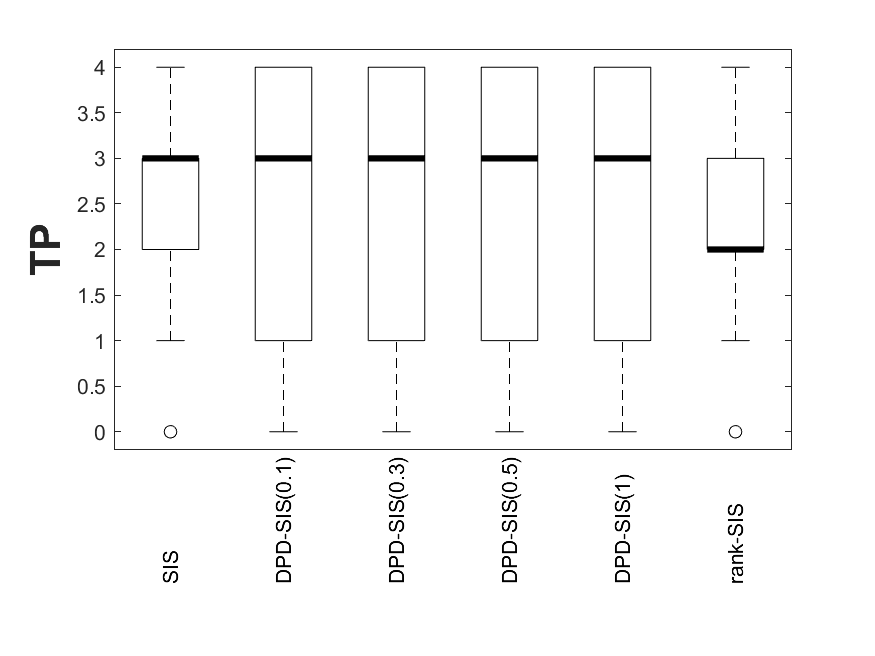}
		\includegraphics[width=0.425\textwidth]{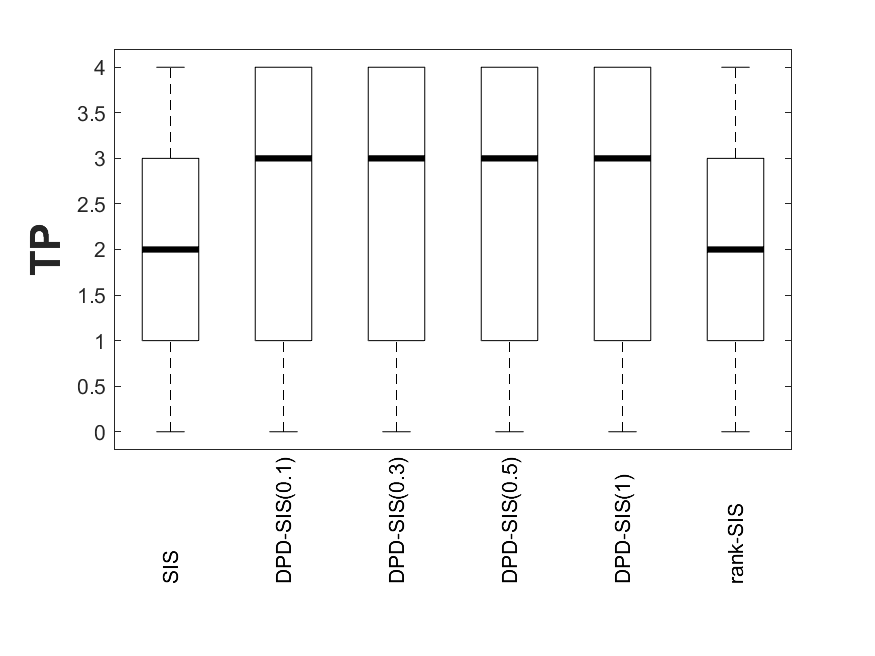}
		\label{FIG:TP_linear_CSIS}}
	\caption{{ The box-plots of true-positives selected by the DPD-CSIS at different $\alpha$ for simulations
		with dependent covariates ($\rho = 0.3$) and strong signal  (non-zero values of $\boldsymbol{\beta}_{0}$ are all 5)
		under pure data (left panel) and 10\% contaminated data (right panel). 
		The results for usual CSIS of \cite{Barut/etc:2016} and rank-SIS of \cite{Li/etc:2012a} are also presented in the same plots for comparisons.} }
	\label{FIG:TP_CSIS}
\end{figure}

{

\subsection{Simulation Results: Performance of DPD-CSIS}

We have also repeated our simulation exercises, as mentioned in Section \ref{SEC:SimSettings},
with some modifications to illustrate the performance of the proposed DPD-CSIS. 
In particular, for this purpose, we assume that the first four covariates are known to be important (each having true regression coefficient 5)
which is used as $\mathcal{X}_\mathcal{C}$; the remaining $p-4$ covariates are used as $\mathcal{X}_\mathcal{D}$ from where the variable screening is performed. 
Among $\mathcal{X}_\mathcal{D}$, $s=4$ sparsely distributed covariates are again assumed to be truly significant  as described in Section \ref{SEC:SimSettings};
the remaining settings of the simulation study are also assumed to be the same as before.
Since conditional screening is more important for dependent covariates, 
we only present the results corresponding to the case $\rho=0.3$ (and strong signal)
in Figures \ref{FIG:TP_logistic_CSIS} and \ref{FIG:TP_linear_CSIS}, respectively, 
for the logistic and the linear regression models under both pure and contaminated data.

We can again observe the advantages of the DPD-CSIS over the usual CSIS \citep{Barut/etc:2016} in terms of robustness, 
producing better variable screening results under data contamination. 
Compared to the rank-SIS, the proposed DPD-CSIS also provides better trade-offs between the robustness and efficiency,
in terms of variable screening under contaminated and pure data set-ups, for both linear and logistic regression models.

\subsection{On the Choices of Tuning Parameters $\alpha$ and $\gamma_n$}
\label{SEC:tuning}


The proposed variable screening procedures, both DPD-SIS and DPD-CSIS, involve two tuning parameters. 
The first one, namely the robustness parameter $\alpha$, controls the trade-off between robustness and efficiency of the marginal MDPDEs 
\citep{Basu/etc:2011,Ghosh/Basu:2013,Ghosh/Basu:2016}. 
As a consequence, as argued in Section \ref{SEC:IF},
the robustness of the proposed DPD-SIS and DPD-CSIS increases with increasing $\alpha>0$. 
However, since the efficiency of the marginal MDPDEs decreases as $\alpha>0$ increases, 
it affects the DPD-SIS or DPD-CSIS with greater variations in the number of true positives for larger $\alpha>0$, 
as seen from our extensive simulation exercises. 
So, we need to choose a proper $\alpha$ to get the optimal trade-off between the performances of DPD-SIS or DPD-CSIS 
under pure and contaminated data (since, in practice, the amount of contamination in a given datasets is often unknown).

As a guiding principle, we can suggest some optimal values of $\alpha$ (depending on the used model and correlation among variables) based on empirical investigations via extensive simulation studies. 
Note that, in our simulations, the performance of DPD-SIS for the cases with independent covariates appears to be almost the same at any $\alpha$ in the range $[0.1, 0.5]$
both for linear and logistic regression models; so any value of $\alpha$ within this range should work well in practice.  
In case of dependent covariates though, their performances can vary significantly over $\alpha$;
an $\alpha$ value around $0.1$ seems to perform the best
(refer to Figures \ref{FIG:TP_logistic_C}, \ref{FIG:TP_lrm_C} and \ref{FIG:TP_CSIS}). 
Hence, as an empirical suggestion, we recommend to apply the proposed DPD-SIS or DPD-CSIS with $\alpha=0.1$ in any practical application 
with linear or logistic regression models involving small to moderate contamination (or no contamination at all).
However, for datasets with higher contamination proportion, we might use a slightly larger $\alpha$ value (e.g., 0.3)
following our theoretical robustness study presented in Section \ref{SEC:IF}. 

Alternatively, one might prefer a data-driven algorithm to chose the optimal $\alpha$ for a given dataset.
For classical regression set-ups, there are a few existing such algorithms for choosing optimal $\alpha$ in the computation of the MDPDEs; 
most notably, the one proposed in \cite{Warwick/Jones:2005} and its recent extension by \cite{Basak/etc:2021}.
These algorithms have also been investigated in the context of linear regression and more general GLMs by \cite{Ghosh/Basu:2013,Ghosh/Basu:2016} and \cite{Basak/etc:2021}.
In the present context of DPD-SIS or DPD-CSIS, we can directly use either of these algorithms for data-driven selection of $\alpha$ 
while computing the marginal MDPDEs and use the resulting optimal MDPDEs to perform variable screening.
However, a major problem with this approach is the possibility to get different optimal $\alpha$ for different marginal MDPDEs 
(associated with different covariates) within the same dataset, leading to inconsistency in the whole variable screening procedure. 
Although we can bypass this obstacle by considering a summary measure from the pool of optimal $\alpha$ values obtained for different covariates, 
this exercise would be extremely time consuming for datasets with larger dimensions and would defeat the whole purpose of variable screening.
So, we recommend to use the empirical suggestions for $\alpha$ while using DPD-SIS (or, DPD-CSIS) for initial (robust) variable screening.
As an alternative, one might also consider the union of the variables selected by the DPD-SIS (or DPD-CSIS) with different possible $\alpha$ values 
in the neighborhood of the empirical suggestion (to ensure that nothing important is missed).

The second tuning parameter $\gamma_n$ controls the amount of false positives in the screening procedure.
The theoretically optimal rate of $\gamma_n$ for controlling the false positives has been seen to be $n^{-2\kappa}$
for some $\kappa>0$ as mentioned in Theorems \ref{THM:DPD_SIS_Sample} and \ref{THM:DPD_SIS_Sample-Cond}, respectively, for the DPD-SIS and the DPD-CSIS.
Note that, this is the exact same rate as the one derived in \cite{Fan/Song:2010} for the usual SIS;
so any existing (data-driven) rule for the selection of $\gamma_n$ for the usual SIS can also be applied for  our DPD-SIS or DPD-CSIS at any $\alpha\geq 0$.
As for the choice of the parameter $\alpha$, one could wish for a data-driven approach to the choice of $\gamma_n$. 
Some efforts have been made lately to establish practical procedures to control the amount of false positives in variable screening procedures, most of them based on data splitting. 
\cite{Guo/etc:2022} presented a general approach to the problem that could be applied also to our DPD-SIS and DPD-CSIS procedures. 
However, in practice it is common to retain a fixed number of predictors, e.g. $[n/\log(n)]$ or $(n-1)$. 
Our preferred approach is to use such a hard thresholding rule followed by a regularized regression procedure where false positives can be controlled by, e.g., 
stability selection \citep{Meinshausen/Buhlmann:2010}. 

}

\section{An Application: NOWAC Lung Cancer Data}\label{SEC:Data_examples}

We will apply our DPD-SIS method to a variable selection problem related to the investigation of potential biomarkers for lung cancer. 
In the Norwegian Women and Cancer (NOWAC) study we have data on 125 lung cancer cases of which 97 had developed metastasis at the time of diagnosis
\citep{Lund/etc:2008,Sandanger/etc:2018}. 
For all these women, we have measures of mRNA in blood some time before diagnosis 
(ranging from 0.3 to 7.9 years before diagnosis, with a median time equal to 4.2 years). 
The goal is to relate the mRNA measurements to the classification of metastatic vs. non-metastatic cancer cases. 
The mRNA measurements are based on the microarray technology, and we have data from a total of 11610 probes. 
Thus, we need to perform some sort of variable selection before running our favorite regression model. Our analysis strategy is as follows: 
First, we run our DPD-SIS procedure ({ with $\alpha=0.1$}) and select the top $n-1$ probes. 
Next, we run a standard logistic regression model with elastic net penalty 
(with the mixing parameter fixed at 0.7 and lambda selected by cross-validation) 
followed by stability selection \citep{Meinshausen/Buhlmann:2010} on the selected probes to reach the final set. 
As a comparison, we run the standard SIS followed by the same standard elastic net logistic regression with stability selection. 
We compare the initial set of $n-1$ probes and the final set after stability selection.

Of the 124 initially selected probes, there was an overlap of { 115}. 
After stability selection with a cut-off at a selection probability equal to 0.7 we were left with { 15} probes based on the robust screening 
and 14 probes based on the ordinary SIS screening, which are reported in Table \ref{TAB:results}. 
We ran logistic regression with elastic net penalty on these 15 vs. 14 probes and calculated the area under the ROC curve (AUC). 
The AUC values were { 0.998} for the 15 robustly selected probes and 0.983 for the 14 probes selected based on ordinary SIS. 
Of course, these values are clearly over-optimistic due to overfitting, 
but they indicate that the predictions based on the robustly selected probes are no worse than those based on the non-robust selection procedure.

\begin{table}[h]
	\caption{The estimated regression coefficients in the final models obtained after stability selection 
		while using the proposed DPD-SIS($\alpha=0.1$) versus the usual SIS. 
		The stability selection probability for each probe is given in the parenthesis.}
	\centering
\begin{tabular}{l|l|r r|r r}\hline
Probe ID & Gene &   \multicolumn{2}{c|}{ DPD-SIS  ($\alpha=0.1$)} & \multicolumn{2}{c}{Usual SIS}\\	\hline\hline
x8tTnl6f115xQSV1X4	&	FER1L5	&	-0.4692	&	(0.96)	&	-0.4311	&	(0.98)	\\
l7tROJgVSRulLNNJ18	&	SERHL	&	-0.5217	&	(0.95)	&	-0.5045	&	(0.92)	\\
Hul6v6J0kV5qD3PA3U	&	TMEM105	&	-0.4952	&	(0.90)	&	-0.4631	&	(0.87)	\\
9Jd97nm3QdLiQXEuzE	&	INVS	&	0.4029	&	(0.93)	&	0.3805	&	(0.86)	\\
BWoVQF900KRd2rRTx8	&	NA	&	0.5487	&	(0.87)	&	0.6041	&	(0.84)	\\
r15YSrezLzSP97mOnU	&	SCARNA14	&	-0.3628	&	(0.78)	&	-0.4609	&	(0.84)	\\
fnRCld.v151SEwlQqk	&	ANO7	&	0.5622	&	(0.89)	&	0.5588	&	(0.79)	\\
rqKhAKFSCyCDqjouuI	&	TBRG1	&	-0.3032	&	(0.87)	&	-0.3044	&	(0.78)	\\
lteivVSVR8WbDmUMBQ	&	HBZ	&	--	&	(0.54)	&	-0.0209	&	(0.78)	\\
Bizm6pN4NAz3jzFe3s	&	SNORD113-7	&	0.4117	&	(0.80)	&	0.4136	&	(0.74)	\\
01\_iIjuJFCnq.s7rqo	&	FLVCR1	&	-0.2828	&	(0.72)	&	-0.3071	&	(0.74)	\\
T9XUjiuZBS6556G.gA	&	MGC23270	&	--	&	(0.67)	&	-0.3512	&	(0.73)	\\
i3urBCqi6u9.0l\_cB0	&	FAM86D	&	0.3580	&	(0.74)	&	0.3468	&	(0.70)	\\
KA5bqKgpHd7QnntF8U	&	NA	&	-0.3219	&	(0.70)	&	-0.2756	&	(0.70)	\\
EEpIrACOOgruuHCpRo	&	CCT6A	&	-0.2740	&	(0.81)	&	--	&	(0.69)	\\
lldNIQXQNde\_jgnE0k	&	ABCB4	&	-0.7045	&	(0.90)	&	--	&		\\
cKfyV6FSSuKdp\_jz3k	&	MAPK7	&	-0.2208	&	(0.74)	&	--	&		\\
		\hline
	\end{tabular}
	\label{TAB:results}
\end{table}

{ 
	Twelve probes were present on both lists; thus, we have five non-overlapping probes worthy of further investigation (see Table \ref{TAB:results}).

The three probes not selected by the usual SIS-procedure, but selected by our DPD-SIS, are interesting. 
The probe EEpIrACOOgruuHCpRo is related to the CCT6A gene, which has been shown to be linked to metastatic non-small cell lung cancer \citep{Zhang/etc:2020}, 
which is the target of our analysis. Furthermore, the two remaining probes (lldNIQXQNde\_jgnE0k and cKfyV6FSSuKdp\_jz3k) are linked to the ABCB4 and the MAPK7 genes. 
Both of these are linked to several types of cancer including lung cancer \citep{Kiehl/etc:2014,  Gavine/etc:2015}.
Thus, all three probes seem highly relevant for our example, and they are obviously important to capture although the usual SIS fails to select them.
On the other hand, two probes were selected by usual SIS but not by our DPD-SIS. 
They were lteivVSVR8WbDmUMBQ connected to the HBZ gene, which has been linked to leukemia but no other cancers, 
and T9XUjiuZBS6556G.gA connected to the MGC23270 gene that has not to our knowledge been linked to cancer at all.

}

This example clearly illustrate that we may miss important variables, probes in our example, by using the standard SIS procedure 
due to the effect of outliers in the data, and this can sometimes lead to spurious findings. 
Our proposed robust DPD-SIS can successfully find the most important probes,   
ignoring any contamination effects present in the data, also without any significant loss in the predictive performance.

\section{Discussion}\label{SEC:conclusion}

This paper presents a robust extension of the usual SIS and the conditional SIS for variable screening under ultra-high dimensional GLMs
using the robust minimum DPD estimators of the marginal regression parameters. 
We show that the proposed DPD-SIS and its conditional version (DPD-SSIS) both enjoy the sure screening property
under a reasonable set of assumptions, in line of what is required by usual SIS. 
More importantly, the proposed DPD-SIS is extremely fast with computational times comparable to those of the usual SIS,
in addition to being robust against data contaminations.
It can be noted that similar robust extensions of SIS can be constructed by other low-dimensional robust estimation procedures for GLMs,
but the validity of the sure screening properties is not guaranteed  under the same rather simple assumptions and with such a low computational time. 
A real data application further reveals the  advantages of our proposed DPD-SIS over the usual SIS 
in identifying the truly important variables from noisy data with no loss in predictive power.

The choice of variable selection criteria is controversial and represents a difficult step in the analysis of high-dimensional data, 
in particular different omics data. In many cases, it is common to reduce the dimensionality of the data by filtering 
the variables before performing any analysis. Even when performing LASSO procedures, 
the data is often filtered beforehand and only a portion of the variables is considered. 
Several approaches are used for variable filtering, most commonly based on the variance or range of the variables, 
or alternatively based on the association with the outcome of interest, as well as more ad-hoc criteria. 
The chosen filtering criteria can clearly impact the analysis, and the results and conclusions might be extremely sensitive to this choice. 
Therefore, it is especially important to base the variable filtering on a clear, objective method that can be reproduced and justified. 
When using SIS based procedures, and DPD-SIS in particular, all variables are evaluated and selected based on their marginal regression coefficient. 
Therefore, DPD-SIS provides a non-arbitrary, reproducible choice of the variables to be included in subsequent analysis, 
and can prove extremely useful in practical applications to, for example, omics data.  
Additionally, the conditional extension provided in this work can prove particularly useful in the context of omics data, 
where it might be relevant to condition on other clinically relevant (non-omics) variables.

We will argue that robust methods, and in particular our robust DPD-based method, are particularly useful in variable filtering situations to avoid the influence of outliers, as traditional checks for outliers become infeasible due to the ultra-high dimensionality.

To control for false discoveries, it is recommendable to perform stability selection in such settings. 
As we illustrate in the applied example, it is quite straightforward to include a stability selection procedure in the DPD-SIS. 
The efficient implementation of DPD-SIS makes this possible without intense computational efforts.

However, a major problem may arise in applications of DPD-SIS when covariates are highly dependent, 
just like for most other (one-step) screening methods.
This problem can be solved by applying the proposed DPD-SIS (or DPD-CSIS) iteratively by removing, in each step, 
the variables selected in the previous steps of the iteration, going through the same philosophy as the iterative SIS \citep{Fan/Lv:2008}.
The performance of such an iterative  DPD-SIS  in successfully addressing the dependent covariate problem has been numerically illustrated 
for the case of linear regression in \cite{Ghosh/Thoresen:2020}; the same can also be investigated for the logistic or other GLMs in future.

Finally, due to the simplicity of the proposed DPD-SIS along with its excellent robust performance, 
it would be natural to extend it to variable screening under more general learning models
like  mixed models or censored regression models with ultra high-dimensional covariates. 
This may be the focus of our future work.


\begin{appendix}

{ 	
\section{Validation of Assumptions for Common GLMs}
\label{APP:Assumption} 

\subsection{Logistic regression}

Let us first consider the logistic regression model, which is a prominent member of the GLM class 
for binary responses having Bernoulli distributions. 
In our notation, for the logistic regression, we have $b(\theta) = \log(1+e^\theta)$ and $g$ is the logit link function,
so that $\mu=E[Y|\theta] = b'(\theta) =  e^\theta /(1+e^\theta)$. 
Therefore, for the computation of the marginal MDPDEs, we have $\xi_\alpha(\theta)=\frac{e^\theta(e^{\alpha\theta}-1)}{(1+e^\theta)^{2+\alpha}}$, 
and hence, the function $\psi_\alpha$ has the form 
$$
\psi_\alpha(y,\theta) = (y - b'(\theta)) \frac{e^{\alpha\theta y}}{(1+e^\theta)^\alpha} - \frac{e^\theta(e^{\alpha\theta}-1)}{(1+e^\theta)^{2+\alpha}}.
$$

Now, in the theoretical analyses of the DPD-SIS, we have first defined the function $B_\alpha$, 
which can be simplified for the present case as

$$
B_\alpha(v(\boldsymbol{x})) = b'(\boldsymbol{x}^T\boldsymbol{\beta}_0)  
- \frac{(e^{\boldsymbol{x}^T\boldsymbol{\beta}_0} - e^{v(\boldsymbol{x})})(e^{\alpha v(\boldsymbol{x})} + e^{v(\boldsymbol{x})})}{
	(1+e^{\boldsymbol{x}^T\boldsymbol{\beta}_0})(1+e^{v(\boldsymbol{x})})^{2+\alpha}}.
$$
Then, one can easily check that $B_\alpha'(\cdot)$ is bounded so that Assumption (B1) of Theorem \ref{THM:DPD_SIS_Pop2} holds
and we have the theoretical justification for the DPD-SIS as a variable screening procedure under the logistic regression model. 
The same can also be verified easily for the DPD-CSIS.

Now, we verify Assumptions (A1)--(A7) to justify the sample-level sure screening property and the false-discovery control of the proposed DPD-SIS and DPD-CSIS.
Recall that Assumptions (A5)--(A7) are exactly the same as those assumed in \cite{Fan/Song:2010} for usual SIS (independent of $\alpha$),
and hence, they hold under mild sufficient conditions as described in the literature of the usual SIS. 
Among other assumptions specific to our DPD-SIS (or DPD-CSIS), we first note that 
Assumption (A1) holds for any $\alpha\geq 0$  with $L_\alpha=1$, the bound for the Bernoulli density, and
Assumptions (A2)-(A3) hold, from the standard literature on the MDPDE under logistic regression \citep{Ghosh/Basu:2016,Basu/etc:2017},
whenever each covariate under screening is independent of the conditioning variables (which is only the intercept in DPD-SIS).
Finally, Assumption (A4) requires a large constant $K_n$ in general, and so we need to carefully choose it in order for
the probability bounds in  Theorems \ref{THM:DPD_SIS_Sample} and \ref{THM:DPD_SIS_Sample-Cond} to make sense. 
In the present case, $b'$ and $B_\alpha$ are bounded and so we can take $k_n$ to be an appropriate finite number (independent of $n$).
So, proceeding as in \cite{Fan/Song:2010}, the optimal order for the sequences in Theorems \ref{THM:DPD_SIS_Sample} and \ref{THM:DPD_SIS_Sample-Cond}
are given by 
$$
K_n = n^{(1-2\kappa)/(\tau+2)}, ~~~~~
R_n = \exp\left( - c_4 n^{\tau(1-2\kappa)/(\tau+2)}\right),
$$
and hence, the (upper) probability bounds in our theorems would make sense as long as  $\log p = o\left(n^{\tau(1-2\kappa)/(\tau+2)}\right)$, 
covering the ultra-high dimensional case with non-polynomial dimensionality as desired.

\subsection{Linear regression with Normal Error}

We can also validate all the assumptions in our theoretical derivations for the usual linear regression model (with standard normal error distribution)
as a special case of our GLM set-up, for which we have $b(\theta) = \theta^2/2$ and $g$ as the identity (link) function. 
This leads to $\xi_\alpha(\theta)=0$, and hence, the $\psi_\alpha$ function in the computation of the marginal MDPDEs has the simplified form 
$$
\psi_\alpha(y,\theta) = \frac{(y - \theta)}{(2\pi)^{\alpha/2}} e^{-\frac{\alpha}{2}(y-\theta)^2}.
$$
Therefore, for this special case, we get 
$$
B_\alpha(v(\boldsymbol{x})) = \boldsymbol{x}^T\boldsymbol{\beta}_0 
- \frac{({\boldsymbol{x}^T\boldsymbol{\beta}_0} - {v(\boldsymbol{x})})}{\sqrt{1+\alpha}(2\pi)^{\alpha/2}}.
$$
Thus, $B_\alpha'$ is a constant and so bounded; 
this implies Assumption (B1) and Theorems \ref{THM:DPD_SIS_Pop1} and \ref{THM:DPD_SIS_Pop2} hold, 
justifying the DPD-SIS for the standard linear regression model at population level. 
Similar justifications for the DPD-CSIS is easy to observe for the present example.

Next, considering the $\alpha$-specific assumptions required for the sample level properties, 
it is straightforward to observe that (A1) holds for any $L_\alpha > (2\pi)^{-\alpha/2}$ 
and Assumptions (A2)--(A3) hold whenever each screening covariate is independent of the conditioning variables.
Finally, to get an idea about the choice of $K_n$ in Assumption (A4), 
keeping the probability bounds in  Theorems \ref{THM:DPD_SIS_Sample} and \ref{THM:DPD_SIS_Sample-Cond} meaningful, 
let us consider $k_n^{(\alpha)} = (1+\alpha)L_\alpha \left[ B(K_n+1)  + K_n^\tau/2m_3\right]$, 
Then, as in \cite{Fan/Song:2010}, assuming $A=\max\{ \tau+4, 3\tau + 2\}$, 
we get the optimal order for the sequences in Theorems \ref{THM:DPD_SIS_Sample} and \ref{THM:DPD_SIS_Sample-Cond} as given by 
$$
K_n = n^{(1-2\kappa)/A}, ~~~~~
R_n = \exp\left( - c_4 n^{\tau(1-2\kappa)/A}\right).
$$
With these choices, the (upper) probability bounds in our Theorems \ref{THM:DPD_SIS_Sample} and \ref{THM:DPD_SIS_Sample-Cond}
would make sense whenever $\log p = o\left(n^{\tau(1-2\kappa)/A}\right)$.
Thus, our theoretical results on the sure screening property and false-discovery control of the proposed DPD-SIS and DPD-CSIS
indeed hold for the example of linear regression models with ultra-high (non-polynomial) dimensional covariates.

\subsection{Poisson regression}

As our final illustration, let us consider the Poisson regression model,
a special case of GLMs corresponding to count responses having Poisson distribution and the log link function. 
In this case, the $\psi_\alpha$ function, for the computation of the marginal MDPDEs, has the form 
$$
\psi_\alpha(y,\theta) = (y - e^\theta) \frac{e^{\alpha \theta y}}{(y!)^\alpha}e^{-\alpha e^\theta} - \xi_\alpha(\theta),
$$
where $\xi_\alpha(\theta)$ has no closed form expression (an infinite sum). 
As a result, the deduced function $B_\alpha(v(\boldsymbol{x}))$ also does not have a closed form expression. 
But, one can check that $b'(\theta) = e^\theta$, and hence, $G_\alpha$ is of the order of an exponential function. 
So, as with the usual SIS in \cite{Fan/Song:2010}, Assumption (B2) in Theorem \ref{THM:DPD_SIS_Pop2} holds
whenever the covariates are bounded or have a light tail (e.g., for sub-Gaussian covariates). 
For such cases, we then have the population level justification for the proposed DPD-SIS as a variable screening procedure 
under the Poisson regression model.

Assumptions (A1)--(A3) specific to the DPD-SIS (or, DPD-CSIS) can also be shown to hold for the Poisson regression case,
as in the previous cases of linear and logistic regression models, whenever the covariates under screening are independent of the conditioning variables and the intercept.
However, to get a reasonable choice of $K_n$ in Assumption (A4), we have to make further assumptions on the covariates. 
For example, if the covariates are bounded, then both $k_n$ and $K_n$ can be taken as finite constants,
and hence, we get $R_n = \exp\left( - c_4 n^{1-2\kappa}\right)$ 
[This is true for any GLM with bounded covariates -- not only the Poisson regression]. 
Then, the (upper) probability bounds in our Theorems \ref{THM:DPD_SIS_Sample} and \ref{THM:DPD_SIS_Sample-Cond} make sense for the ultra-high dimensional case
with $\log p = o\left(n^{1-2\kappa}\right)$ and the sample-level properties of the proposed DPD-SIS and DPD-CSIS also hold for such Poisson regression models. 
We conjecture that the same can be verified for sub-Gaussian covariates as well (as observed empirically) 
although we do not have a concrete proof at this moment. 

In this respect, we must note that, there is indeed no literature on the validity of such an assumption, even for the usual SIS,  under Poisson regression with more general covariate distributions.
So, more research is needed to establish the validity of the required assumptions for SIS (and also for our proposed DPD-SIS or DPD-CSIS)  under such general Poisson models.

}

\section{Proofs of Theorems \ref{THM:DPD_SIS_Pop1} and \ref{THM:DPD_SIS_Pop2}}

\subsection{Proof of Theorem \ref{THM:DPD_SIS_Pop1}}
	This follows from the Fisher consistency of the MDPDEs (by definition)
	and the results of Theorem 2 in \cite{Fan/Song:2010}.

	\subsection{Proof of Theorem \ref{THM:DPD_SIS_Pop2}}

	For $\alpha=0$ this theorem is identical to Theorem 3 of \cite{Fan/Song:2010}.
	We now extend their argument to prove it for any given $\alpha>0$.
	Let us fix an $\alpha>0$ and $j \in \mathcal{M}_0$.
	
	First consider the cases where $B_\alpha'(\cdot)$ is bounded and let $D$ be its upper bound.
	Then $B_\alpha(\cdot)$ is Lipschitz continuous and hence 
	$$
	\left|\left\{B_\alpha(\beta_{j0}^{M\alpha}+\beta_j^{M\alpha}X_j) - B_\alpha(\beta_{j0}^{M\alpha})\right\}X_j\right|
	\leq D|\beta_j^{M\alpha}|X_j^2.
	$$ 
	Taking expectation, we get 
	\begin{eqnarray}
		D|\beta_j^{M\alpha}| &\geq& \left|E\left[\left\{B_\alpha(\beta_{j0}^{M\alpha}+\beta_j^{M\alpha}X_j) - B_\alpha(\beta_{j0}^{M\alpha})\right\}X_j\right]\right|
		\nonumber\\
		&=& \left|E\left[B_\alpha(\beta_{j0}^{M\alpha}+\beta_j^{M\alpha}X_j)X_j\right]\right|
		= \left|\mbox{Cov}\left(B_\alpha(\beta_{j0}^{M\alpha}+\beta_j^{M\alpha}X_j), X_j\right)\right|.
		\nonumber
	\end{eqnarray}
{ In the above, we have used $E\left[B_\alpha(\beta_{j0}^{M\alpha})X_j\right]=0$
which holds since $B_\alpha(\beta_{j0}^{M\alpha})$ is a constant (by definition) and $EX_j=0$ by our assumption.}
	But, from the estimating equation in (8), we get 
	$
	E[B_\alpha(\beta_{j0}^{M\alpha}+\beta_j^{M\alpha}X_j) X_j]  = E[b'(\boldsymbol{X}^T\boldsymbol{\beta}_0)X_j],
	$
	and hence, using $E(X_j)=0$, 
	\begin{eqnarray}
		\mbox{Cov}(B_\alpha(\beta_{j0}^{M\alpha} +\beta_j^{M\alpha}X_j), X_j)  = \mbox{Cov}(b'(\boldsymbol{X}^T\boldsymbol{\beta}_0), X_j).
		\label{EQ:app2}
	\end{eqnarray}
	Then, in view of the condition of the theorem, we get $|\beta_j^{M\alpha}| \geq D_1^{-1}c_1 n^{-\kappa}$
	completing the proof of the theorem.
	
	Next, we consider the cases where the second assumption, namely Condition (9), holds.
	Clearly if $|\beta_j^{M\alpha}| \geq c n^{-\kappa}$ for a sufficiently large universal constant $c>0$, 
	the result holds and we are done. So, assume $|\beta_j^{M\alpha}| \leq \widetilde{c}_1 n^{-\kappa}$ for some $\widetilde{c}>0$
	and let $\beta_0^{M\alpha}$ be a constant such that $B_\alpha(\beta_0^{M\alpha}) = E[Y]$.
	We first prove the following claim.
	
	\noindent
	Claim 1: $\left|\beta_{j0}^{M\alpha} - \beta_0^{M\alpha}\right|\leq \widetilde{c}_2$ for all $j\in \mathcal{M}_0$ and some constant $\widetilde{c}_2>0$.
	
	To prove the claim, we fix a $j\in\mathcal{M}_0$ and consider the marginal MDPDE objective function (population version) as a function of $\beta_0$ only as
	$Q(\beta_0) = E[l_\alpha(Y, \beta_0 + \beta_j^{M\alpha}X_j)]$ so that we get 
	$$
	Q'(\beta_0) = E[Y - B_\alpha(\beta_0 + \beta_j^{M\alpha}X_j)] = B_\alpha(\beta_0^{M\alpha}) - E[B_\alpha(\beta_0 + \beta_j^{M\alpha}X_j)].
	$$
	But, 
	\begin{eqnarray}
		&& \left|E[B_\alpha(\beta_0 + \beta_j^{M\alpha}X_j)] - B_\alpha(\beta_0)\right|
		\nonumber\\
		&\leq& \sup_{|x|\leq \widetilde{c}_1n^{\eta - \kappa}}\left|B_\alpha(\beta_0 + x) - B_\alpha(\beta_0)\right|
		+ 2E\left[G_\alpha(a|X_j|)|X_j|I(|X_j|> n^{\eta})\right]
		\nonumber\\
		&=& o(1) + o(1),\nonumber
	\end{eqnarray}
	by the continuity of $B_\alpha(\cdot)$ and Condition (9).
	Therefore, we get $Q'(\beta_0) = B_\alpha(\beta_0^{M\alpha}) - B_\alpha(\beta_0) + o(1)$
	and hence, for a $\widetilde{c}_2>0$, we have $Q'(\beta_0^{M\alpha}-\widetilde{c}_2) <0$ and 
	$Q'(\beta_0^{M\alpha}+\widetilde{c}_2) >0$ since $B_\alpha(\cdot)$ is strictly increasing.
	Hence  $\left|\beta_{j0}^{M\alpha} - \beta_0^{M\alpha}\right|\leq \widetilde{c}_2$ proving our Claim 1.

	Finally, to prove the theorem, we note that if $|X_j|\leq n^\kappa$, then Claim 1 ensures that the points
	$\beta_{j0}^{M\alpha}$ and $(\beta_{j0}^{M\alpha}+ \beta_j^{M\alpha}X_j)$, for all $j\in \mathcal{M}_0$,  
	belong the interval $I=(\beta_0^{M\alpha}-h, \beta_0^{M\alpha}+h)$ independent of  $j$, where $h=\widetilde{c}-1 + \widetilde{c}_2$. 
	Let $\widetilde{D}=\max_{x\in I} B_\alpha'(x)$, which is finite by Lipschitz continuity of $B_\alpha(\cdot)$
	in a neighborhood of $\beta_0^{M\alpha}$ and hence, for $|X_j|\leq n^\kappa$, we have
	$$
	\left|\left\{B_\alpha(\beta_{j0}^{M\alpha}+ \beta_j^{M\alpha}X_j) - B_\alpha(\beta_{j0}^{M\alpha})\right\}\right|
	\leq \widetilde{D}|\beta_j^{M\alpha}|X_j^2.
	$$
	Taking expectation over the region $\{|X_j|\leq n^\kappa\}$, we get
	\begin{eqnarray}
		\widetilde{D}|\beta_j^{M\alpha}| &\geq& 
		\left|E\left[\left\{B_\alpha(\beta_{j0}^{M\alpha}+\beta_j^{M\alpha}X_j) - B_\alpha(\beta_{j0}^{M\alpha})\right\}X_jI(|X_j|\leq n^\kappa)\right]\right|
		\nonumber\\
		&=& \left|\mbox{Cov}(b'(\boldsymbol{X}^T\boldsymbol{\beta}_0), X_j)\right| - A_0 - A_1,
		\label{EQ:app3}
	\end{eqnarray}
	by a similar calculation leading to (\ref{EQ:app2}), where 
	$A_m = E\left[B_\alpha(\beta_{j0}^{M\alpha}+\beta_j^{M\alpha}X_j^m)X_jI(|X_j|> n^\kappa)\right]$ for $m=0,1$.
	But, $\left|\beta_{j0}^{M\alpha}+\beta_j^{M\alpha}X_j^m\right|\leq a|X_j|$ for $|X_j|>n^\kappa$ with a sufficiently large $n$ independent of $j$ and $m$,
	we get from Condition (9)  that 
	$
	A_m \leq E[G(a|X_j|)^m|X_j|I(|X_j|\geq n^\kappa)] \leq dn^{-\kappa},
	$
	for both $m=01,1$. Then, the theorem follows from (\ref{EQ:app3}) using the given condition that 
	$\left|\mbox{Cov}(b'(\boldsymbol{X}^T\boldsymbol{\beta}_0), X_j)\right|\geq c_1n^{-\kappa}$. 
	\hfill{$\square$}

\section{Proof of Lemma \ref{LEM:MDPDE_expCov}}
	
	The result in the lemma holds directly by Theorem 1 of \cite{Fan/Song:2010}, provided we can show that their Conditions (A), (B) and (C) 
	are implied by our Assumptions (A1)--(A5). 
	In this regard, note that Assumption (A3) is indeed a reformulation of Condition  (A) of \cite{Fan/Song:2010}.
	Further, Assumption (A2) implies Condition (C) of \cite{Fan/Song:2010} via a second order Taylor series expansion of 
	$l_\alpha\left(Y, \boldsymbol{X}_{j}^T\boldsymbol{\beta}_j\right)$ with respect to $\boldsymbol{\beta}_j$ 
	around $\boldsymbol{\beta}_j = \boldsymbol{\beta}_j^{M\alpha}$. Finally it remains to show that Condition (B) of \cite{Fan/Song:2010}
	holds under Assumptions (A1), (A4) and (A5).

	Let us define $\Omega_n = \left\{ (X_j, Y): |X_j|\leq K_n, |Y| \leq K_n^* \right\}$, where $K_n$ is as in Assumption (A4) and 
	$K_n^\ast = \frac{m_0}{m_3}K_n^\tau$ with $m_0, m_3$ and $\tau$ being as in Assumption (A5). Then, for our present case, 
	Condition (B) of \cite{Fan/Song:2010} becomes equivalent to \\
	\begin{eqnarray}
		\left|l_\alpha(Y, \boldsymbol{X}_{j}^T\boldsymbol{\beta}_j) - l_\alpha(Y, \boldsymbol{X}_{j}^T\boldsymbol{\beta}_j')\right|I((X_j,Y)\in\Omega_n)
		\leq k_n^{(\alpha)} \left|\boldsymbol{X}_{j}^T\boldsymbol{\beta}_j - \boldsymbol{X}_{j}^T\boldsymbol{\beta}_j'\right|I((X_j,Y)\in\Omega_n),
		\nonumber\\
		~~~~~\boldsymbol{\beta}_j, \boldsymbol{\beta}_j'\in\mathcal{B},
		\label{EQ:Cond_B1}
		\\
		\mbox{and }~~
		\sup\limits_{\boldsymbol{\beta}_J\in\mathcal{B}: ||\boldsymbol{\beta}_j-\boldsymbol{\beta}_j^{M\alpha}||\leq \epsilon_1}
		\left|E\left[l_\alpha(Y, \boldsymbol{X}_{j}^T\boldsymbol{\beta}_j) - l_\alpha(Y, \boldsymbol{X}_{j}^T\boldsymbol{\beta}_j^{M\alpha})\right]
		I((X_j,Y)\notin\Omega_n)\right|\leq o(n^{-1}),~~
		\label{EQ:Cond_B2}
	\end{eqnarray}
	where $k_n^{(\alpha)}$ is as defined in the statement of the Lemma and $\epsilon_1$ as in Assumption (A4).
	First, to show (\ref{EQ:Cond_B1}), we use a first order Taylor series expansion to get
	\begin{eqnarray}
		l_\alpha(Y, \boldsymbol{X}_{j}^T\boldsymbol{\beta}_j) - l_\alpha(Y, \boldsymbol{X}_{j}^T\boldsymbol{\beta}_j')
		= D(\widetilde{\boldsymbol{\beta}}_j) 
		\left[\boldsymbol{X}_{j}^T\boldsymbol{\beta}_j - \boldsymbol{X}_{j}^T\boldsymbol{\beta}_j'\right],
		\label{EQ:Lemma_P1}
	\end{eqnarray}
	where $\widetilde{\boldsymbol{\beta}}_j\in\mathcal{B}$ lies on the line segment joining $\boldsymbol{\beta}_j$ and  $\boldsymbol{\beta}_j'$
	and \\
	$
	D(\widetilde{\boldsymbol{\beta}}_j) =(1+\alpha)\left[\xi_\alpha(\boldsymbol{X}_j^T\widetilde{\boldsymbol{\beta}}_j)
	- (Y- b'(\boldsymbol{X}_j^T\widetilde{\boldsymbol{\beta}}_j))f^\alpha(Y;\boldsymbol{X}_j^T\widetilde{\boldsymbol{\beta}}_j)\right].
	$
	But, on $\Omega_n$, we have
	\begin{eqnarray}
		\left|D(\widetilde{\boldsymbol{\beta}}_j)\right| 
		&\leq& (1+\alpha)\left[\left|\xi_\alpha(\boldsymbol{X}_j^T\widetilde{\boldsymbol{\beta}}_j)\right|
		+ (|Y|+ |b'(\boldsymbol{X}_j^T\widetilde{\boldsymbol{\beta}}_j)|)\left|f^\alpha(Y;\boldsymbol{X}_j^T\widetilde{\boldsymbol{\beta}}_j)\right|\right]
		\nonumber\\
		&\leq& (1+\alpha)\left[\left|\xi_\alpha(K_nB+B)\right|+ \left(\frac{m_0}{m_3}K_n^\tau + |b'(K_nB+B)|\right)L_\alpha\right]=k_n^{(\alpha)},
		\nonumber
	\end{eqnarray}
	by Assumption (A1), (A5) and the subsequent result in (12). 
	Substituting it in (\ref{EQ:Lemma_P1}), we get Condition (\ref{EQ:Cond_B1}).
	
	Next, to prove (\ref{EQ:Cond_B2}), we again consider the expansion (\ref{EQ:Lemma_P1}) with $\boldsymbol{\beta}_j'= \boldsymbol{\beta}_j^{M\alpha}$
	and by taking expectation we get
	\begin{eqnarray}
		&&\left|E\left[l_\alpha(Y, \boldsymbol{X}_{j}^T\boldsymbol{\beta}_j) - l_\alpha(Y, \boldsymbol{X}_{j}^T\boldsymbol{\beta}_j^{M\alpha})\right]\right|
		= E\left|D(\widetilde{\boldsymbol{\beta}}_j) \left[\boldsymbol{X}_{j}^T\boldsymbol{\beta}_j - \boldsymbol{X}_{j}^T\boldsymbol{\beta}_j^{M\alpha}\right]\right|
		\nonumber\\
		&&\leq(1+\alpha)||\boldsymbol{\beta}_j - \boldsymbol{\beta}_j^{M\alpha}||_2 E\left|\left[B_\alpha(\boldsymbol{X}_{j}^T\widetilde{\boldsymbol{\beta}}_j) - B_\alpha(\boldsymbol{X}_{j}^T\boldsymbol{\beta}_j^{M\alpha})\right]||\boldsymbol{X}_j||_2\right|,
		\nonumber
	\end{eqnarray}
	by an application of the Cauchy-Schwartz inequality. Therefore,  we get
	\begin{eqnarray}
		&& \sup\limits_{\boldsymbol{\beta}_j\in\mathcal{B}: ||\boldsymbol{\beta}_j-\boldsymbol{\beta}_j^{M\alpha}||\leq \epsilon_1}
		\left|E\left[l_\alpha(Y, \boldsymbol{X}_{j}^T\boldsymbol{\beta}_j) - l_\alpha(Y, \boldsymbol{X}_{j}^T\boldsymbol{\beta}_j^{M\alpha})\right]
		I((X_j,Y)\notin\Omega_n)\right|
		\nonumber\\
		&&
		\leq(1+\alpha)\epsilon_1 \sup\limits_{\boldsymbol{\beta}_j\in\mathcal{B}: ||\boldsymbol{\beta}_j-\boldsymbol{\beta}_j^{M\alpha}||\leq \epsilon_1}
		E\left[|B_\alpha(\boldsymbol{X}_{j}^T\widetilde{\boldsymbol{\beta}}_j)|||\boldsymbol{X}_j||_2 
		+ |B_\alpha(\boldsymbol{X}_{j}^T\boldsymbol{\beta}_j^{M\alpha})|||\boldsymbol{X}_j||_2\right]I(|X_j|>K_n),
		\nonumber\\
		&& \leq o(n^{-1}),\nonumber
	\end{eqnarray}
	by Assumption (A4), and this completes the proof.
	\hfill{$\square$}

\section{Proof of Theorem \ref{THM:DPD_SIS_Sample}}
	
	\noindent\textbf{Part (a):}\\
	We start with Lemma \ref{LEM:MDPDE_expCov} and take $(1+t)= c_3Vn^{\frac{1}{2}-\kappa} (16k_n^{(\alpha)})^{-1}>0$ to get
	\begin{eqnarray}
		P\left(\left|\widehat{\beta}_j^{M\alpha} - \beta_j^{M\alpha}\right|\geq c_3n^{-\kappa} \right)
		\leq e^{-\frac{n^{1-2\kappa}}{K_n^2K_n^2}C} + nm_1e^{-m_0K_n^\tau} = R_n,
		~~~~~j=1, \ldots, p,
		\label{EQ:Pf_t1}
	\end{eqnarray}
	Then, the uniform convergence result in Part (a) of the theorem holds from the relation (\ref{EQ:Pf_t1})
	via  union bound of probabilities. 
	
	\bigskip
	\noindent\textbf{Part (b):}\\
	Let us consider the event 
	$\mathcal{E}_n = \left\{ \max\limits_{j \in\mathcal{M}_0}\left|\widehat{\beta}_j^{M\alpha} - \beta_j^{M\alpha}\right| \leq c_2n^{-\kappa}/2 \right\}$.\\
	By Theorem \ref{THM:DPD_SIS_Pop2}, on $\mathcal{E}_n$, we then have $\left|\widehat{\beta}_j^{M\alpha}\right|\geq c_2n^{-\kappa}/2$ for all $j\in\mathcal{M}_0$.
	Therefore, for the choice of $\gamma_n$ as given in the statement of the theorem, we have $\mathcal{M}_0 \subset \widehat{\mathcal{M}}_\alpha(\gamma_n)$
	on $\mathcal{E}_n$, and hence 
	$$
	P\left(\widehat{\mathcal{M}}(\gamma_n) \supset \mathcal{M}_0 \right)\geq P(\mathcal{E}_n) = 1 - P(\mathcal{E}_n^c).
	$$
	But, since $\mathcal{M}_0$ has $s$ elements, by a union bound of probability, we get from (\ref{EQ:Pf_t1}) that 
	$P(\mathcal{E}_n^c) \leq sR_n$ completing the proof of Part (b).

	\bigskip
	\noindent\textbf{Part (c):}\\
	The proof is based on the result (13), evaluated at $\boldsymbol{\beta}_j=\boldsymbol{\beta}_j^{M\alpha}$, 
	which implies that the number of variables having $|\beta_j^{M\alpha}|>\epsilon n^{-\kappa}$ cannot exceed $O(n^{2\kappa}\Lambda_{\max}(\Sigma))$
	for any given $\epsilon>0$. Now, let us consider the event
	$$
	\widetilde{\mathcal{E}}_n = \left\{\max\limits_{1\leq j \leq p} |\widehat{\beta}_j^{M\alpha} - \beta_j^{M\alpha}| \leq \epsilon n^{-\kappa} \right\}.
	$$ 
	Then, on the event $\widetilde{\mathcal{E}}_n$, we have 
	$$
	\left| \left\{ j :|\widehat{\beta}_j^{M\alpha}|>2\epsilon n^{-\kappa}  \right\}\right|
	\leq \left| \left\{ j :|{\beta}_j^{M\alpha}|>2\epsilon n^{-\kappa}  \right\}\right| \leq O(n^{2\kappa}\Lambda_{\max}(\Sigma)).
	$$
	Hence, taking $\epsilon=c_5/2$ for the choice of $\gamma_n$ as given in the statement of the theorem, we get 
	$
	P\left(|\widehat{\mathcal{M}}(\gamma_n)| \leq O(n^{2\kappa}\Lambda_{\max}(\Sigma)) \right)\geq P(\widetilde{\mathcal{E}}_n)
	=1-P(\widetilde{\mathcal{E}}_n^c).
	$
	But, by Part (a) of the theorem, we have $P(\widetilde{\mathcal{E}}_n^c)\leq pR_n$ completing the proof.
	\hfill{$\square$}

\section{Proofs of Theorems \ref{THM:DPD_SIS_Pop1_Cond} and \ref{THM:DPD_SIS_Pop2_Cond}}
	
We note that, for each $j\in\mathcal{D}$, the quantity $\boldsymbol{\beta}_{\mathcal{C}j}^{M\alpha}$, 	defined in (17), 
	satisfies the estimating equations given by 
	\begin{eqnarray}
		E \left[\psi_\alpha\left(Y, \boldsymbol{X}_{\mathcal{C}j}^T\boldsymbol{\beta}_{\mathcal{C}j}^{M\alpha}\right)\boldsymbol{X}_{\mathcal{C}}\right] 
		= 0,
		~~~~~~~~
		E \left[\psi_\alpha\left(Y, \boldsymbol{X}_{\mathcal{C}j}^T\boldsymbol{\beta}_{\mathcal{C}j}^{M\alpha} \right)X_j\right] = 0.
		\label{EQ:Est_eqn_Mpop_Cond}
	\end{eqnarray}
	On the other hand, the baseline quantity $\boldsymbol{\beta}_{\mathcal{C}}^{M\alpha}$ satisfies 
	\begin{eqnarray}
		E \left[\psi_\alpha\left(Y, \boldsymbol{X}_{\mathcal{C}}^T\boldsymbol{\beta}_{\mathcal{C}}^{M\alpha}\right)\boldsymbol{X}_{\mathcal{C}}\right] 
		= 0.
		\label{EQ:Est_eqn_Mpop_Cond0}
	\end{eqnarray}
	Further, for any $j\in\mathcal{D}$,  using $E[X_j|\boldsymbol{X}_{\mathcal{C}}]=0$, we have 
	\begin{eqnarray}
		\mbox{Cov}_L(Y, X_j|\boldsymbol{X}_{\mathcal{C}})
		&=& E\left[\left(Y-E[Y|\boldsymbol{X}_{\mathcal{C}}]\right)X_j\right]
		= E\left[YX_j\right] 
		\nonumber\\
		&=& E\left[E(Y|\boldsymbol{X})X_j\right] 
		=E \left[b'(\boldsymbol{X}^T\boldsymbol{\beta}_0)X_j\right],
		\label{EQ:cov_L0}
	\end{eqnarray}
	and hence, invoking the definitions of $\psi_\alpha$ and $B_\alpha$,  we get
	\begin{eqnarray}
		E \left[\psi_\alpha\left(Y, \boldsymbol{X}_{\mathcal{C}}^T\boldsymbol{\beta}_{\mathcal{C}}^{M\alpha} \right)X_j\right] 
		&=& E\left[\left(B_\alpha(\boldsymbol{X}_{\mathcal{C}}^T\boldsymbol{\beta}_{\mathcal{C}}^{M\alpha}) - b'(\boldsymbol{X}^T\boldsymbol{\beta}_0)\right)
		X_j\right]
		\nonumber\\
		&=& E E\left[\left(B_\alpha(\boldsymbol{X}_{\mathcal{C}}^T\boldsymbol{\beta}_{\mathcal{C}}^{M\alpha}) - b'(\boldsymbol{X}^T\boldsymbol{\beta}_0)\right)
		X_j|\boldsymbol{X}_{\mathcal{C}}\right]
		\nonumber\\
		&=& - E \left[b'(\boldsymbol{X}^T\boldsymbol{\beta}_0)X_j\right]
		\nonumber\\
		&=& - \mbox{Cov}_L(Y, X_j|\boldsymbol{X}_{\mathcal{C}}).
		\label{EQ:Est_eqn_Mpop_Cond4}
	\end{eqnarray}

\subsection{Proof of Theorem \ref{THM:DPD_SIS_Pop1_Cond}}
	
	Firstly, if $\beta_j^{M\alpha}=0$ for some $j\in\mathcal{D}$, from (\ref{EQ:Est_eqn_Mpop_Cond}) we get
	\begin{eqnarray}
		E \left[\psi_\alpha\left(Y, \boldsymbol{X}_{\mathcal{C}}^T\boldsymbol{\beta}_{\mathcal{C}j1}^{M\alpha}\right)\boldsymbol{X}_{\mathcal{C}}\right] 
		= 0,
		~~~~~~~~
		E \left[\psi_\alpha\left(Y, \boldsymbol{X}_{\mathcal{C}}^T\boldsymbol{\beta}_{\mathcal{C}j1}^{M\alpha} \right)X_j\right] = 0.
		\label{EQ:Est_eqn_Mpop_Cond2}
	\end{eqnarray}
	Combining the first equation with (\ref{EQ:Est_eqn_Mpop_Cond0}) and the uniqueness of its solution we have 
	$\boldsymbol{\beta}_{\mathcal{C}j1}^{M\alpha}  = \boldsymbol{\beta}_{\mathcal{C}}^{M\alpha} $
	and hence the second equation in (\ref{EQ:Est_eqn_Mpop_Cond2}) becomes
	\begin{eqnarray}
		E \left[\psi_\alpha\left(Y, \boldsymbol{X}_{\mathcal{C}}^T\boldsymbol{\beta}_{\mathcal{C}}^{M\alpha} \right)X_j\right] = 0.
		\label{EQ:Est_eqn_Mpop_Cond3}
	\end{eqnarray}
	This leads to the desired condition $\mbox{Cov}_L(Y, X_j|\boldsymbol{X}_{\mathcal{C}})=0$ by (\ref{EQ:Est_eqn_Mpop_Cond4}).

	On the other hand, if  $\mbox{Cov}_L(Y, X_j|\boldsymbol{X}_{\mathcal{C}})=0$ for some $j\in\mathcal{D}$, 
	then by (\ref{EQ:Est_eqn_Mpop_Cond4}), the equation in (\ref{EQ:Est_eqn_Mpop_Cond3}) hold.
	Combining (\ref{EQ:Est_eqn_Mpop_Cond3}) with (\ref{EQ:Est_eqn_Mpop_Cond0}), we see that 
	$\boldsymbol{\beta}_{\mathcal{C}j}^{M\alpha} = \left(\boldsymbol{\beta}_{\mathcal{C}}^{M\alpha}, 0\right)^T$ is a solution 
	of the estimating equations in (\ref{EQ:Est_eqn_Mpop_Cond}), leading to $\beta_j^{M\alpha}=0$.
	\hfill{$\square$}

\subsection{Proof of Theorem \ref{THM:DPD_SIS_Pop2_Cond}}

	Fix any $j\in\mathcal{M}_{0\mathcal{D}}$ and 
	define $\boldsymbol{\Omega}_j = E\left[m_{\alpha,j} \boldsymbol{X}_{\mathcal{C}j}\boldsymbol{X}_{\mathcal{C}j}^T\right]$
	and $\boldsymbol{\beta}_{\Delta, j}= \left(\boldsymbol{\beta}_{\mathcal{C}j1}^{M\alpha} - \boldsymbol{\beta}_{\mathcal{C}}^{M\alpha}\right)$.
	Consider a partition of $\boldsymbol{\Omega}_j$ as given by 
	$$
	\boldsymbol{\Omega}_j = \begin{bmatrix}
		\begin{array}{cc}
			\boldsymbol{\Omega}_{11,j} & \boldsymbol{\Omega}_{12,j}\\
			\boldsymbol{\Omega}_{21,j}^T & \boldsymbol{\Omega}_{22,j}
		\end{array}
	\end{bmatrix}
	= \begin{pmatrix}
		\begin{array}{cc}
			E\left[m_{\alpha,j} \boldsymbol{X}_{\mathcal{C}}\boldsymbol{X}_{\mathcal{C}}^T\right] & E\left[m_{\alpha,j} \boldsymbol{X}_{\mathcal{C}}{X}_j\right]
			\\
			E\left[m_{\alpha,j} {X}_{j}\boldsymbol{X}_{\mathcal{C}}^T\right] & E\left[m_{\alpha,j} {X}_{j}^2\right]
		\end{array}
	\end{pmatrix}.
	$$
	
	Now, from the estimating equations (\ref{EQ:Est_eqn_Mpop_Cond}) and (\ref{EQ:Est_eqn_Mpop_Cond0}), 
	along with the definitions of $B_\alpha$ and $m_{\alpha, j}$, we get 
	\begin{eqnarray}
		0 &=& E\left[\left(B_\alpha(\boldsymbol{X}_{\mathcal{C}j}^T\boldsymbol{\beta}_{\mathcal{C}j}^{M\alpha}) 
		- B_\alpha(\boldsymbol{X}_{\mathcal{C}}^T\boldsymbol{\beta}_{\mathcal{C}}^{M\alpha})\right)\boldsymbol{X}_{\mathcal{C}}\right]
		\nonumber\\
		&=& E\left[m_{\alpha,j}\left(\boldsymbol{X}_{\mathcal{C}j}^T\boldsymbol{\beta}_{\mathcal{C}j}^{M\alpha} - 
		\boldsymbol{X}_{\mathcal{C}}^T\boldsymbol{\beta}_{\mathcal{C}}^{M\alpha}\right)\boldsymbol{X}_\mathcal{C}\right]
		\nonumber\\
		&=& E\left[m_{\alpha,j}\left(\boldsymbol{X}_{\mathcal{C}}^T\boldsymbol{\beta}_{\Delta,j} + {X}_{j}{\beta}_{j}^{M\alpha}\right)
		\boldsymbol{X}_\mathcal{C}\right].
		\nonumber
	\end{eqnarray}
	Therefore, by solving, we get
	$\boldsymbol{\beta}_{\Delta,j} = - \Omega_{11,j}^{1}\Omega_{12,j}\beta_j^{M\alpha}$.
	Further, note that, for any integrable function $h(\boldsymbol{X}_\mathcal{C})$, 
	we have 
	$$
	E[h(\boldsymbol{X}_\mathcal{C})X_j] = EE[h(\boldsymbol{X}_\mathcal{C})X_j|\boldsymbol{X}_\mathcal{C}]
	= E\left[h(\boldsymbol{X}_\mathcal{C})E(X_j|\boldsymbol{X}_\mathcal{C})\right] = 0.
	$$
	Hence, by (\ref{EQ:cov_L0}), (\ref{EQ:Est_eqn_Mpop_Cond}), and the definition of $m_{\alpha,j}$, we get 
	\begin{eqnarray}
		\mbox{Cov}_L(Y, X_j|\boldsymbol{X}_{\mathcal{C}})
		&=& E\left[b'(\boldsymbol{X}^T\boldsymbol{\beta}_0)X_j\right]  = E\left[\left(b'(\boldsymbol{X}^T\boldsymbol{\beta}_0)
		-B_\alpha(\boldsymbol{X}_\mathcal{C}^T\boldsymbol{\beta}_{\mathcal{C}}^{M\alpha})\right)X_j\right]
		\nonumber\\
		&=& E\left[\left(B_\alpha(\boldsymbol{X}_{\mathcal{C}j}^T\boldsymbol{\beta}_{\mathcal{C}j}^{M\alpha})
		-B_\alpha(\boldsymbol{X}_\mathcal{C}^T\boldsymbol{\beta}_{\mathcal{C}}^{M\alpha})\right)X_j\right]
		\nonumber\\
		&=& E\left[m_{\alpha,j}\left(\boldsymbol{X}_{\mathcal{C}j}^T\boldsymbol{\beta}_{\mathcal{C}j}^{M\alpha}
		- \boldsymbol{X}_\mathcal{C}^T\boldsymbol{\beta}_{\mathcal{C}}^{M\alpha}\right)X_j\right]
		\nonumber\\
		&=& \Omega_{12,j}^T\boldsymbol{\beta}_{\Delta,j} + \Omega_{22,j}\beta_j^{M\alpha}
		\nonumber\\
		&=& \left[\Omega_{22,j} - \Omega_{12,j}^T\Omega_{11,j}^{-1}\Omega_{12,j}\right]\beta_j^{M\alpha}.
		\nonumber
	\end{eqnarray}
	Now, taking absolute value in the above and using the assumptions of the theorem, we get
	$
	c_1n^{-\kappa} \leq c_2 \left|\beta_j^{M\alpha}\right|.
	$
	Since this holds for all $j\in\mathcal{M}_{0\mathcal{D}}$, taking minimum over all such $j$ 
	we get the desired conclusion of the theorem with $c_3 = c_1/c_2$.
	\hfill{$\square$}
\end{appendix}

\section*{Acknowledgment}
This research is partially supported by the Norwegian Research Council - grant numbers 248804 and 262111. 
The NOWAC post-genome cohort study was funded by the ERC advanced grant Transcriptomic in Cancer Epidemiology (ERC-2008-AdG-232997).
A major part of this research work has been done while the first author (AG) was visiting University of Oslo, Norway. 
The research of AG is also supported by an INSPIRE Faculty research grant and a grant (No. SRG/2020/000072) from SERB,
both under the Department of Science and Technology, Government of India


%


\begin{thebibliography}{}
	
	
	
	\bibitem[\protect\citeauthoryear{Barut }{Barut et al.}{2016}]{Barut/etc:2016}
	Barut, E., Fan, J., and Verhasselt, A. (2016). 
	Conditional sure independence screening. 
	\textit{J Amer Stat Assoc}, 111(515), 1266-1277.
	
	\bibitem[\protect\citeauthoryear{Basak..}{Basak et~al.}{2021}]{Basak/etc:2021}
	Basak, S., Basu, A., and Jones, M. C. (2021). 
	On the `optimal' density power divergence tuning parameter. 
	\textit{J Appl Stat}, 48(3), 536--556.
	
	
	\bibitem[\protect\citeauthoryear{Basu, Harris, Hjort and Jones}{Basu et~al.}{1998}]{Basu/etc:1998}
	Basu, A., Harris, I. R., Hjort, N. L., and Jones, M. C. (1998). 
	Robust and efficient estimation by minimising a density power divergence. 
	{\it Biometrika}, {\bf 85}, 549--559.
	
	
	\bibitem[\protect\citeauthoryear{Basu, Shioya and Park}{Basu et~al.}{2011}]{Basu/etc:2011}
	Basu, A., Shioya, H. and Park, C. (2011). 
	\emph{Statistical Inference: The Minimum Distance Approach}. 
	Chapman $\&$ Hall/CRC, Boca de Raton.
	
	
	
	\bibitem[\protect\citeauthoryear{Basu, Mandal, Martin and Pardo}{Basu et~al.}{2017}]{Basu/etc:2017}
	Basu, A., Ghosh, A., Mandal, A., Martin, N. and Pardo, L. (2017)
	A Wald-type test statistic for testing linear Hypothesis in logistic regression models based on minimum density power divergence estimator.
	{\it Electron. J. Stat.}, {\bf 11}, 2741--2772.
	
	
	\bibitem[\protect\citeauthoryear{Basu, Mandal, Martin and Pardo}{Basu et~al.}{2021}]{Basu/etc:2018}
	Basu, A., Ghosh, A., Mandal, A., Martin, N., and Pardo, L. (2021). 
	Robust Wald-type tests in GLM with random design based on minimum density power divergence estimators. 
	\textit{Stat Method Appl}, 30(3), 973--1005.
	
	
	
	\bibitem[\protect\citeauthoryear{Buhlmann, P., and Van De Geer}{Buhlmann and Van De Geer}{2011}]{Buhlmann/VanDeGeer:2011}
	Buhlmann, P., and Van De Geer, S. (2011). 
	\textit{Statistics for high-dimensional data: methods, theory and applications}. 
	Springer Science \& Business Media.
	
	
	
	\bibitem[\protect\citeauthoryear{Fan and Li}{Fan and Li}{2001}]{Fan/Li:2001}
	Fan, J. and Li, R. (2001). 
	Variable Selection via Nonconcave Penalized Likelihood and its Oracle Properties. 
	\textit{J Amer Statist Assoc}, 96:1348–1360.
	
	
	\bibitem[\protect\citeauthoryear{Fan and Lv}{Fan and Lv}{2008}]{Fan/Lv:2008}
	Fan, J., and Lv, J. (2008). Sure independence screening for ultrahigh dimensional feature space. 
	\textit{J Royal Stat Soc B},, 70(5), 849-911.
	
	\bibitem[\protect\citeauthoryear{Fan and Song}{Fan and Song}{2010}]{Fan/Song:2010}
	Fan, J., and Song, R. (2010). 
	Sure independence screening in generalized linear models with NP-dimensionality. 
	\textit{Ann Stat}, 38(6), 3567-3604.
	
	
	\bibitem[\protect\citeauthoryear{Fu, L., and Wang}{Fu and Wang}{2018}]{Fu/Wang:2018}
	Fu, L., and Wang, Y. G. (2018). 
	Variable selection in rank regression for analyzing longitudinal data. 
	\textit{Stat Methods Med Res}, 27(8), 2447-2458.
	
	
	\bibitem[\protect\citeauthoryear{Gather, U., and Guddat}{Gather and Guddat}{2008}]{Gather/Guddat:2008}
	Gather, U., and Guddat, C. (2008). 
	Comment on “Sure Independence Screening for Ultrahigh Dimensional Feature Space” by Fan, JQ and Lv, J. 
	\textit{J Royal Stat Soc B}, 70, 893-895.
	
	{
	
	\bibitem[\protect\citeauthoryear{Gavine...}{Gavine et al.}{2015}]{Gavine/etc:2015}
	Gavine, P. R., Wang, M., Yu, D., Hu, E., Huang, C., Xia, J., ... and Ji, Q. (2015). 
	Identification and validation of dysregulated MAPK7 (ERK5) as a novel oncogenic target in squamous cell lung and esophageal carcinoma. 
	\textit{BMC cancer}, 15(1), 1--9.

}
	
	\bibitem[\protect\citeauthoryear{Ghosh }{Ghosh}{2019}]{Ghosh:2019}
	Ghosh, A. (2019). 
	Robust inference under the beta regression model with application to health care studies. 
	\textit{Stat. Method. Med Res}, 28(3), 871--888.
	
	
	\bibitem[\protect\citeauthoryear{Ghosh and Basu}{Ghosh and Basu}{2013}]{Ghosh/Basu:2013}
	Ghosh, A., and Basu, A. (2013). 
	Robust estimation for independent non-homogeneous observations using density power divergence with applications to linear regression. 
	{\em Electron. J. Stat.}, {7}, 2420--2456. 
	
	
	%
	
	\bibitem[\protect\citeauthoryear{Ghosh and Basu}{Ghosh and Basu}{2016}]{Ghosh/Basu:2016}
	Ghosh, A., and Basu, A. (2016). 
	Robust Estimation in Generalized Linear Models : The Density Power Divergence Approach. 
	{\em Test}, {25(2)}, 269--290.
	
	
	
	
	
	%
	
	
	
	\bibitem[\protect\citeauthoryear{Ghosh and Majumdar}{Ghosh and Majumdar}{2020}]{Ghosh/Majumdar:2019}
	Ghosh, A. and Majumdar, S. (2020). 
	Ultrahigh-dimensional Robust and Efficient Sparse Regression using Non-Concave Penalized Density Power Divergence. 
	\textit{IEEE Trans. Info. Theory}, 66(12), 7812--7827.
	
	\bibitem[\protect\citeauthoryear{Ghosh and Thoresen}{Ghosh and Thoresen}{2021}]{Ghosh/Thoresen:2020}
	Ghosh, A. and Thoresen, M. (2021). 
	A Robust Variable Screening procedure for Ultra-high dimensional data. 
	\textit{Stat. Meth. Med. Res.},  30(8), 1816--1832.
	
	
	
	\bibitem[\protect\citeauthoryear{Giraud}{Giraud}{2014}]{Giraud:2014}
	Giraud, C. (2014). \textit{Introduction to high-dimensional statistics}. 
	Chapman and Hall/CRC.
	
{ 
	\bibitem[\protect\citeauthoryear{Guo...}{Guo et al.}{2022}]{Guo/etc:2022}
	Guo, X., Ren, H., Zou, C. and Li R. (2022).
	Threshold selection in feature screening for error rate control.
	\textit{J Amer Statist Assoc}, doi: 10.1080/01621459.2021.2011735.
}
	
	
	\bibitem[\protect\citeauthoryear{Hall, P., and Miller}{Hall and Miller}{2009}]{Hall/Miller:2009}
	Hall, P., and Miller, H. (2009). 
	Using generalized correlation to effect variable selection in very high dimensional problems. 
	\textit{J Comput Graphical Stat}, 18(3), 533-550.
	
	
	\bibitem[\protect\citeauthoryear{Hampel, Ronchetti, Rousseeuw, and  Stahel}{Hampel et~al.}{1986}]{Hampel/etc:1986}
	Hampel, F.~R., Ronchetti, E., Rousseeuw, P.~J., and Stahel W.(1986).
	\newblock {\em Robust Statistics: The Approach Based on Influence Functions}.
	\newblock New York, USA: John Wiley \& Sons.
	
	
	\bibitem[\protect\citeauthoryear{Hastie}{Hastie et~al.}{2015}]{Hastie/etc:2015}
	Hastie, T., Tibshirani, R., and Wainwright, M. (2015). 
	\textit{Statistical learning with sparsity: the lasso and generalizations}. 
	CRC press.
	
	
	
	
	
	
	\bibitem[\protect\citeauthoryear{Kiehl}{Kiehl et al.}{2014}]{Kiehl/etc:2014}
	Kiehl, S., Herkt, S. C., Richter, A. M., Fuhrmann, L., El-Nikhely, N., Seeger, W., ..., and Dammann, R. H. (2014). 
	ABCB4 is frequently epigenetically silenced in human cancers and inhibits tumor growth. 
	\textit{Sci. Rep.}, 4(1), 1--9.
	
	
	\bibitem[\protect\citeauthoryear{Li }{Li et al.}{2012a}]{Li/etc:2012a}
	Li, G., Peng, H., Zhang, J., and Zhu, L. (2012a). 
	Robust rank correlation based screening. 
	\textit{Ann Stat}, 40(3), 1846-1877.
	
	\bibitem[\protect\citeauthoryear{Li }{Li et al.}{2012b}]{Li/etc:2012b}
	Li, R., Zhong, W., and Zhu, L. (2012b). 
	Feature screening via distance correlation learning. 
	\textit{J Amer Statist Assoc}, 107(499), 1129-1139.
	
	
	\bibitem[\protect\citeauthoryear{Lund et al}{Lund et al.}{2008}]{Lund/etc:2008}
	Lund, E., Dumeaux, V., Braaten, T., Hjart{a}ker, A., Engeset, D., Skeie, G., and Kumle, M. (2008). 
	Cohort profile: the Norwegian women and cancer study—NOWAC—Kvinner og kreft. 
	\textit{Int J Epidem.}, 37(1), 36--41.
	
	
	\bibitem[\protect\citeauthoryear{Luo et al}{Luo et al.}{2014}]{Luo/etc:2014}
	Luo, S., Song, R., and Witten, D. (2014). 
	Sure Screening for Gaussian Graphical Models. 
	\textit{Stat}, 1050, 29.4
	
	\bibitem[\protect\citeauthoryear{Meinshausen and Buhlmann}{Meinshausen and Buhlmann}{2010}]{Meinshausen/Buhlmann:2010}
	Meinshausen, N. and Buhlmann, P. (2010). 
	Stability Selection. 
	\textit{J Royal Stat Soc B}, 72, 417-473.
	
	
	
	\bibitem[\protect\citeauthoryear{Mu, W., and Xiong}{Mu and Xiong}{2014}]{Mu/Xiong:2014}
	Mu, W., and Xiong, S. (2014). 
	Some notes on robust sure independence screening. 
	\textit{J App Stat}, 41(10), 2092--2102.
	
	
	
	
	
	
	
	
	\bibitem[\protect\citeauthoryear{Saldana, D. F., and Feng}{Saldana and Feng}{2018}]{Saldana/Feng:2018}
	Saldana, D. F., and Feng, Y. (2018). 
	SIS: An R package for sure independence screening in ultrahigh-dimensional statistical models. 
	\textit{J Stat Software}, 83(2), 1--25.
	
	
	\bibitem[\protect\citeauthoryear{Lund et al}{Sandanger et al.}{2018}]{Sandanger/etc:2018}
	Sandanger, T. M., N\o st, T. H., Guida, F., Rylander, C., Campanella, G., Muller, D. C., ... and Chadeau-Hyam, M. (2018). 
	DNA methylation and associated gene expression in blood prior to lung cancer diagnosis in the Norwegian Women and Cancer cohort. 
	\textit{Sci. Rep.}, 8(1), 1--10.
	
	
	\bibitem[\protect\citeauthoryear{van der Vaart}{van der Vaart}{1998}]{vanderVaart}
	van der Vaart, A. W. (1998). 
	\newblock {\em Asymptotic statistics}. 
	Cambridge University Press.
	
	\bibitem[\protect\citeauthoryear{Wainwright}{Wainwright}{2019}]{Wainwright:2019}
	Wainwright, M. J. (2019). 
	\textit{High-dimensional statistics: A non-asymptotic viewpoint} (Vol. 48). Cambridge University Press.
	
	\bibitem[\protect\citeauthoryear{Wang}{Wang et al.}{2017}]{Wang/etc:2017}
	Wang, T., Zheng, L., Li, Z., and Liu, H. (2017). 
	A robust variable screening method for high-dimensional data. 
	\textit{J App Stat}, 44(10), 1839-1855.
	
	\bibitem[\protect\citeauthoryear{Warwick and Jones}{Warwick and	Jones}{2005}]{Warwick/Jones:2005}
	Warwick, J. and M.~C. Jones (2005).
	\newblock Choosing a robustness tuning parameter.
	\newblock {\em J Stat Comput Simul\/}, {75}, 581--588.
	
	{ 
	
	\bibitem[\protect\citeauthoryear{Zhang...}{Zhang et al.}{2020}]{Zhang/etc:2020}
	Zhang, T., Shi, W., Tian, K., and Kong, Y. (2020). 
	Chaperonin containing t-complex polypeptide 1 subunit 6A correlates with lymph node metastasis, 
	abnormal carcinoembryonic antigen and poor survival profiles in non-small cell lung carcinoma. 
	\textit{World journal of surgical oncology}, 18(1), 1--10.

}

	\bibitem[\protect\citeauthoryear{Zhao/Li}{Zhao and Li}{2012}]{Zhao/Li:2012}
	Zhao, S. D., and Li, Y. (2012). 
	Principled sure independence screening for Cox models with ultra-high-dimensional covariates. 
	\textit{J Mult Anal}, 105(1), 397-411.
	
	
	\bibitem[\protect\citeauthoryear{Zhong}{Zhong}{2014}]{Zhong:2014}
	Zhong, W. (2014). 
	Robust sure independence screening for ultrahigh dimensional non-normal data. 
	\textit{Acta Mathematica Sinica}, English Series, 30(11), 1885--1896.
	
	
	
\end{thebibliography}
\end{document}